
\font\steptwo=cmr10 scaled\magstep2

\magnification=\magstep1
\settabs 18 \columns
\hoffset=.20truein
\hsize=16truecm


\def\Cit{\hbox{\it l\hskip -5.5pt C\/}}

\def\Crm{\hskip0.5mm \hbox{\rm l\hskip -5.pt C\/}}
\def\Zrm{Z\hskip-2mm Z}

\baselineskip=17 pt

\def\s{\smallskip}

\def\b{\bigskip}
\def\bb{\bigskip\bigskip}

\def\no{\noindent}

\def\ce{\centerline}
\def\ve{\vfill\eject}
 
\def\YB{Y\hskip-1mmB}

\font\got=eufm9 scaled\magstep1
\def\g{{\got g}}

\def\harr#1#2{\smash{\mathop{\hbox to .25 in{\rightarrowfill}}
 \limits^{\scriptstyle#1}_{\scriptstyle#2}}}

\def\today{\ifcase\month\or January\or February\or March\or April\or 
May\or June\or July\or August\or September\or October\or November\or 
December\fi \space\number\day, \number\year }




\def\Cit{\hbox{\it l\hskip -5.5pt C\/}}

{\ce {\bf  GENERALIZATION AND EXACT DEFORMATIONS  }}
\b
{\ce{\bf OF QUANTUM GROUPS}}
 
\b
 {\ce {C. Fr\o nsdal}}

{\ce { \it Department of Physics and Astronomy}} 
{\ce {\it University of California, Los Angeles CA 90095-1547, USA}}
\b
 \bb\b
\no {\bf Abstract.}
A large family of ``standard" coboundary Hopf algebras is investigated. 
The existence of a universal R-matrix is demonstrated for the case when
the parameters are in general position. Algebraic surfaces in parameter
space are characterized by the appearance of certain ideals; in this case the
universal R-matrix exists on the associated algebraic quotient. In
special cases the quotient is a ``standard" quantum group; all familiar
quantum groups including twisted ones are obtained in this way. In other
special cases one finds new types of coboundary bi-algebras.

The ``standard" universal R-matrix is shown to be the unique solution of a
very simple, linear recursion relation. The classical limit is obtained in the case of
quantized Kac-Moody algebras of finite and affine type. 

 Returning
to the general case, we study deformations of the standard R-matrix and the 
associated Hopf algebras. A preliminary investigation of the first order deformations
uncovers a class of deformations that incompasses  the quantization of all 
Kac-Moody algebras of finite and affine type. The corresponding exact deformations
are described as generalized twists, $ R_\epsilon =
(F^t)^{-1}RF$, where $R$ is the standard R-matrix and the cocycle $F$ (a power series in 
the deformation parameter $\epsilon$) is
the solution of a linear recursion relation of the same type as that which determines
$R$. Included here is the universal R-matrix for the elliptic quantum groups associated with 
$sl(n)$, a big surprise! 

Specializing again, to
the case of quantized  Kac-Moody algebras, and taking the  classical limit of these
esoteric quantum groups, one re-discovers all the trigonometric and elliptic  r-matrices of Belavin and
Drinfeld.  The formulas obtained here are easier to use than the original ones, and
the structure of the space of classical r-matrices  
is more transparent.  
The r-matrices obtained here are more general in that they are defined on the full Kac-Moody algebras,
the central
extensions of the loop groups.

\ve
\ce {\steptwo {\bf TABLE}}
\b
1. ~Introduction.

2. ~Standard Universal R-matrices.

3. ~Differential Algebras.

4. ~Differential Complexes.

5. ~Integrability of Eq.(2.14).

6. ~Completion of the Proof of Theorem 2.

7. ~Obstructions and Generalized Serre Relations.

8. ~The Standard Classical r-matrix for Simple Lie Algebras.

9. ~The Standard Classical r-matrix for Untwisted Loop Algebras.

10. The Standard Classical r-matrix for Twisted Loop Algebras.

11. Including the Central Extension.
\b

\ce {{\it DEFORMATIONS}}

\b
12. First order Deformations.

13. First order Deformations of Type $e_\sigma \otimes e_{-\rho}$ and the Classical Limit.

14. Hopf Structure.

15. Exact Deformations of Standard, Generalized Quantum Groups.

16. Esoteric r-matrices.

18. Universal Elliptic R- and r-matrices.

~~~~~Acknowledgements.

~~~~~References.

\ve
\voffset0truein
\line {\bf 1. Introduction.   \hfil}

Quantum groups sprouted in that fertile soil where mathematics overlaps
with physics.  The mathematics of quantum groups is exciting, and the
applications to physical modelling are legion.  It is the more surprising
that some aspects of the structure of quantum groups remain to be
explored; this is especially true of those aspects that bear upon the
problem of classification.  The quantum groups that have so far found
employment in physics are very special (characterized by a single
``deformation" parameter
$q$).  It is true that these applications are susceptible to some
generalization, by the process of ``Cartan twisting"; by this we mean the 
type of twisting that was used by Reshetikhin [R] to construct the multiparameter quantum groups, 
in which a quantum R-matrix $R$ is 
replaced by $\tilde R = (F^t)^{-1}RF$, with F in the Cartan subalgebra.  Unfortunately it is easy to
receive the impression that twisting is a gauge
transformation that relates equivalent structures.  The fact that Cartan-twisted
or multiparameter quantum groups differ qualitatively among themselves
becomes evident when one investigates their rigidity to deformation. 
Deformation theory is a means of attacking the classification problem; at
the same time it offers a wider horizon against which to view the whole
subject.  The new quantum groups discovered this way (the deformations of
the twisted ones) are dramatically different; the physical applications
should be of a novel kind.

Let {\g}  be a simple Lie algebra over \Crm.  A structure of
coboundary Lie bialgebra on {\g } is determined by a ``classical"
r-matrix; an element $r\in ~$\g$ \,\otimes \,${\g} that satisfies the
classical Yang-Baxter relation
$$  [r_{12},r_{13} + r_{23}] + [r_{13},r_{23}]=0~, \eqno(1.1)
$$
\no as well as the symmetry condition
$$
 r+r^t=\hat K~, \eqno(1.2)
$$
\no where $\hat K$ is the Killing form of \g.  The classification
of r-matrices of simple complex Lie algebras (finite and affine) was accomplished by Belavin
and Drinfeld [BD].

It is widely believed that there corresponds, to each such r-matrix, via
a process of ``quantization,"  a unique quantum group [D2]. Somewhat more
precisely, one expects that there exists a Hopf algebra deformation
$\tilde U($\g$)$ of $U($\g$)$, and an element
$R\in \tilde U($\g$)\otimes \tilde U($\g$)$, such that $\Delta
R=R\Delta^\prime$, where $\Delta$ is the coproduct of $\tilde U($\g$)$
and
$\Delta^\prime$ is the opposite coproduct,   satisfying the (quantum)
Yang-Baxter relation
$$
 R_{12}R_{13}R_{23}=R_{23}R_{13}R_{12}~; \eqno(1.3)
$$
\no such that $r$ can be recovered by an expansion of $R$ with respect to
a parameter $\hbar$:
$$  
R=1+\hbar r+o(\hbar^2)~. \eqno(1.4)
$$
Till now, this program had been realized for r-matrices of a  
class that we  call ``standard".
\vskip.5cm

\no {\bf Definition 1.1.}  Let \g ~be a simple, complex Lie algebra,
\g$^0$ a Cartan subalgebra and $\Delta^+$ a set of positive roots.
A (constant) standard r-matrix for \g~ has the expression
$$
 r=r_0 + \sum_{\alpha\in\Delta^+} E_{-\alpha}\otimes E_\alpha~. \eqno(1.5)
$$
\no Here $r_0\in $~\g$^0 \,\otimes\, $\g$^0$ is restricted by (1.2). An affine r-matrix 
(non-constant, with spectral parameter)   is of standard type if it commutes with the 
Cartan subalgebra.
\vskip.5cm

The (universal) R-matrix that corresponds to a standard r-matrix is
known. Explicit formulas are of two types: in terms of Serre generators, 
or in terms Lie generators. An explicit formula in terms of Lie generators has been 
given for the simplest choice of
$r_0$ in [KR].  An expression for $R$ in terms of Drinfeld-Serre generators [D1][FR][LS][Ro][T]
seems more fundamental (especially so in the affine case),
$$
 R=R^0\bigl(1+\sum_\alpha e_{-\alpha}\otimes e_\alpha + \ldots\bigr)~.
\eqno(1.6)
$$
\no Here $\{H_a,e_\alpha,e_{-\alpha}\}$ are Chevalley-Drinfeld generators
associated with a Cartan subalgebra and simple roots,
$R^0$ involves only the $H_a$'s.  An R-matrix of this form will be called
standard; a precise   definition (in a more general context) will be
given  in Section \nobreak 2, \penalty-1000   Definition 2.2.  The relationship between (1.5) and
(1.6) is examined in Sections 8 and 16. An explicit formula for the coefficients in (1.6) is in (5.9).  

The R-matrices associated with the multiparameter quantum groups discovered by
Reshetikhin [R] and others [Sc][Su] are thus all included in the rubrique
``standard".  The principal characteristic of a standard R-matrix is that
it ``commutes with Cartan": 
$$
\bigl[H_a\otimes 1+1\otimes H_a,R\bigr]=0.
$$  Until now, non-standard R-matrices were known only in the fundamental
representation [CG][FG1].
\ve

The aim of this work is to use deformation theory to discover
 the so far unknown quantum groups that are alleged
to be associated with non-standard r-matrices.  This seems a reasonable
approach because (i) non-standard r-matrices can be viewed, and
effectively calculated [F], as deformations of standard  r-matrices and
(ii) the largest family of non-standard quantum groups known so far was
found by applying deformation theory to  certain standard R-matrices in
the fundamental representation [FG2].

Progress  achieved 
in the present paper  is due, in the first place,
to the idea of focusing on the representation (1.6) of the standard
universal R-matrix, and in the second place to the discovery of a
differential complex associated with the Yang-Baxter relation: the  
study of (1.6) turned out to be unexpectedly rewarding.

The existence of a universal R-matrix in the form (1.6), 
for quantized Lie algebras and for Kac-Moody algebras, was known [D1][FR][T].  But it
  turns out that the representation (1.6) for an R-matrix that satisfies
the Yang-Baxter equation makes sense in a context that is much wider than
quantized Kac-Moody algebras. We present a proof of the existence of an R-matrix of the form (1.6)
that covers a   wider category of bialgebras. The proof is constructive and
provides  useful insight into the structure of these bialgebras (actually Hopf
algebras). It exploits  a direct connection between the Yang-Baxter equation and
a certain differential complex, and it reduces  the calculation of $R$ to the
solution of a linear recursion relation.

We introduce (Definition 2.1) an algebra ${\cal{A}}$  with generators
$\{H_a,e_\alpha,e_{-\alpha}\}$ that satisfy certain relations, including
the following (see also Eq. (1.7b) below):
$$
 [H_a,H_b]=0,~[H_a,e_{\pm\alpha}] = \pm H_a(\alpha) e_{\pm\alpha}~,\quad H_a(\alpha) \in \Crm.
\eqno(1.7a)
$$
\no    We define a standard R-matrix on
${\cal{A}}$~--~Definition 2.2.~--~as a formal series of the form
$$
\eqalign{R&=\exp(\varphi^{ab}H_a\otimes H_b)\bigl(1+e_{-\alpha}
\otimes e_\alpha  + \sum^\infty_{k=2} t^{(\alpha^\prime)}_{(\alpha)}
e_{-\alpha_1}\ldots e_{-\alpha_k}\otimes e_{\alpha^\prime_1}\ldots
e_{\alpha^\prime_k}\bigr)~, \cr}
$$
\no with parameters $\varphi^{ab}\in \Crm$, fixed, and   determine
the coefficients $t^{(\alpha^\prime)}_{(\alpha)}\in  \Crm
 $ so that the Yang-Baxter relation (1.3) is satisfied.  One finds that
this requires additional relations, namely
$$  [e_\alpha,e_{-\beta}] = \delta^\beta_\alpha
\bigl(e^{\varphi(\alpha,\cdot)}-e^{-\varphi(\cdot,\alpha)}\bigr)~, 
\eqno(1.7b)
$$
\no with $\varphi(\alpha,\cdot)=\varphi^{ab}H_a(\alpha)H_b$,
$\varphi(\cdot,\alpha)=\varphi^{ab}H_aH_b(\alpha)$.  These relations are
therefore included in the definition of the algebra
${\cal{A}}$. Generically, with the parameters in
general position, no further relations are required.

The generators $H_a$ of ${\cal{A}}$ generate an Abelian subalgebra 
  ${\cal{A}}^0$    that may be called the Cartan subalgebra. ~ A
key point is to refrain from introducing,  {\it a priori},   any
(generalized) Serre  relations 
among the  Chevalley-Drinfeld generators
$e_\alpha$, or  among the
$e_{-\alpha}$.  The algebras of ultimate interest are obtained
subsequently, by identifying an appropriate ideal $I\subset {\cal{A}}$
that intersects ${\cal{A}}^0$ trivially, and passing to the quotient ${\cal A}' = 
{\cal{A}}/I$.  This is the strategy of Chevalley [C], fully exploited in
the theory of Kac-Moody algebras [K][Mo]; here it is applied to 
``generalized quantum groups."

This point of view allows a significant generalization. We study free differential algebras
in general, then attempt to classify the ideals. To each ideal there corresponds a quotient 
algebra on which a coboundary Hopf algebra (with its standard R-matrix) can be constructed. 
Quantized Kac-Moody algebras, characterized by  Serre-Drinfeld ideals, form a special case.

The first result is Theorem 2.  It asserts that the Yang-Baxter relation
for the standard R-matrix on $\cal A$ is equivalent to a simple, linear
recursion relation for the coefficients
$t^{(\alpha^\prime)}_{(\alpha)}$. This result is of great help in the subsequent calculations.   

The integrability of this recursion relation, Eq.(2.14), is related to the first 
cohomology group on quantum planes. Generically, all one-forms are exact, whence 
the second result that, when the parameters of ${\cal{A}}$ are in general position, 
there exists a unique set of coefficients
$t^{(\alpha^\prime)}_{(\alpha)}$ such that the standard R-matrix
satisfies the Yang-Baxter relation.

Obstructions to the solution of the recursion
relation (2.14), and thus to the Yang-Baxter relation on ${\cal A}$, exist on certain hyper-surfaces 
in the space of parameters of
${\cal{A}}$; they are detected by the presence of ``constants". A Serre relation is a special 
type of constant. Constants are
   studied in a   general context in Sections 3 and 4; their
complete classification is an open, but probably not unsolvable
problem.  Its solution would shed light on the structure of ideals in quantized Kac-Moody 
algebras and reduce the theorem of Gabber and Kac [GK] to a corollary. 
The relevance of this discussion to the Yang-Baxter relation is demonstrated
in Section 5, and the proof of Theorem 2 is completed in Section 6. 

The study of the obstructions is taken up again in Section 7.  The third
main result is Theorem 7: the obstructions (that is, the constants)
generate an ideal $I\subset {\cal{A}}$, and a unique standard R-matrix,
satisfying Yang-Baxter, exists on ${\cal{A}}/I$.
 
Next,   we specialize to quantized Kac-Moody algebras and calculate the  classical limit,
that is, the classical r-matrix associated with the standard R-matrix.  The result was of
course known, but without a precise determination of $R$ it is not possible  to evaluate the
limit directly. A further complication is that all the coefficients become singular. But the recursion
relation (2.14) guarantees the existence of the limit and provides an efficient method for
evaluating it.  See Sections 9,10 and 16,17.
 
  Quantized Kac-Moody algebras are characterized by the
property that, for each pair $(\alpha,\beta)$, there is a positive integer 
$k = k_{\alpha\beta}$ such that the following relation holds
$$ 
e^{\varphi(\alpha,\beta) + \varphi(\beta,\alpha) + (k-1)\varphi(\alpha,\alpha)}
 =1.\eqno(1.8)
$$ 
In this case the ideal $I$ is generated by the Serre relations
$$ 
0 = \sum^k_{m=0} Q^k_m~(e_\alpha)^me_\beta\, (e_\alpha)^{k-m}~,
\eqno(1.9)
$$ 
with coefficients
 $$
 Q^k_m=(-)^m e^{m\varphi(\alpha,\beta)}~q^{m(m-1)/2}
\left(\matrix{k\cr m\cr}\right)_q~,\quad q := e^{\varphi(\alpha,\alpha)}.  
$$ 

\no We suppose Card $N$ and Card $M$ finite and interpret $A = 1-k$ as the generalized Cartan
matrix of a Kac-Moody algebra.

The classical r-matrix associated with $R$ is defined after a rescaling of the generators  as the 
coefficient of $\hbar$ in  the expansion $R = 1 + \hbar r + o(\hbar^2)$; it satisfies the classical Yang-Baxter relation. 
(Note that $r + r^t \neq 0$; the antisymmetric part of $r$ satisfies the modified 
classical Yang-Baxter relation.) 
We  calculate this classical r-matrix, dealing separately  with the following cases: 
 Kac-Moody algebras of finite type in Section 8;   unextended  loop algebras 
(untwisted and twisted) in Sections 9 and 10;    the  full Kac-Moody
algebras in Section 11.
\ve

\ce {DEFORMATIONS}
\vskip.5cm

The rest of the paper is a study of the deformations of the standard R-matrix, satisfying the Yang-Baxter relation,
but in the wider context of the bialgebras 
${\cal A}$ and ${\cal A}' = {\cal A}/I$ described above. We set 
$$
R_\epsilon = R + \epsilon R_1 + o(\epsilon^2),
$$
and suppose that  $R_1$ is driven by a term of the type
$$
Se_{-\rho} \otimes e_\sigma + S'e_\sigma \otimes e_{-\rho},\,\,\,\,\, S,S' \in {\cal A}'^0.
\eqno(1.10)
$$
Such deformations exist under certain conditions on the parameters; then $S$ and $S'$ and
the remaining terms in
$R_\epsilon$  
  (a formal power series in $\epsilon$ with constant
term $R$) are determined by the Yang-Baxter relation.

  We begin by calculating a class
of first order deformations of the standard R-matrix on ${\cal{A}}/I$ for
any ideal of obstructions
$I\subset {\cal{A}}$.  This is our fourth result, Theorem 13.1. The main
difficulty is that the problem is not well posed, for we have been unable
to discover a category that is both natural and convenient in which to
calculate \underbar {all} deformations.  We  limit our study to a class
of deformations.  The good news is Theorem 13.2: when we specialize   to
the case of simple quantum groups, then we obtain quantizations of all  
simple Lie bialgebras (constant r-matrices) of Belavin and Drinfeld, so far to first order in the deformation parameter.

In Section 14 we define the coproduct, counit and antipode that
turn all these algebras into Hopf algebras. This completes the investigation of first order deformations.
The results provide inspiration for construction of exact deformations.

  An exact formula (to all orders in
$\epsilon$) in closed form  for $R_\epsilon$  is obtained for the case of elementary
deformations, when
$R_1$ is a single term  of the type (1.10). In the general case of  compound
deformations, when (1.10) is replaced by a sum of terms of the same type,  we obtain exact
deformations  in the form of a generalized twist. (Section 15.)

\ve
Let $R$ be the R-matrix of a coboundary Hopf algebra ${\cal A}'$, and $F \in {\cal A}' \otimes {\cal A}'$, invertible.
Then
$$
\tilde R := (F^t)^{-1}RF\eqno(1.11)
$$
satisfies the Yang-Baxter relation if $F$ satisfies the cocycle condition [G]
 
$$
\bigl((1 \otimes \Delta_{21}) F\bigr)F_{12} = \bigl(
(\Delta_{13} \otimes 1) F\bigr) F_{31}.\eqno(1.12)
$$
(See Theorem 15.1 for the complete statement.) Though it is not quite germaine to our discussion, 
it may be worth while to point out that, if $R$ is unitary, then so is $\tilde R$; the formula (1.11)
therefore yields a  family of (mostly) new unitary R-matrices. 

Applying this to our context, we find that the  relation (1.12) is equivalent
to a simple, linear recursion relation that can be reduced to the same form as the 
recursion relation that determines the coefficients in the expansion of $R$.
  It has a unique solution that can be expressed directly in terms of these coefficients.
   Just as in the standard case, this leads to a
simple equation for the  classical r-matrix, from which the  latter is determined to all
orders.  

In Section 16 we specialize to the case of deformed quantized, affine Kac-Moody algebras   and  take
the classical limit, to recover the esoteric, affine r-matrices of the simple Lie algebras, with 
their central extensions. The result agrees with that of Belavin and Drinfeld, except that they 
did not include the central extension. The formulas obtained in this paper are
more transparent and simpler to use.

Finally, in Section 17, we deal with a very special case, to discover that the elliptic r-matrices of 
$sl(N)$ also arise as the classical limit of certain deformed quantum groups.
 The universal elliptic R-matrix is 
expressed as an infinite product. It is shown, in the particular case of the elliptic R-matrix for 
$sl(2)$ in the fundamental representation, that this infinite product is both convergent and 
of practical utility; it reduces to the representation of elliptic functions 
in terms of infinite products, and the result is in perfect agreement with Baxter [B].

\ve

\line {{\bf 2. Standard Universal R-matrices.} \hfil}
\s 

The universal R-matrix of a standard or twisted quantum group has the form
$$
\eqalign{R&=\exp(\varphi^{ab}H_a\otimes H_b)  \cr &\,\,\,\,\times
\bigl(1+t_\alpha(e_{-\alpha}\otimes e_\alpha)+t_{\alpha\beta}
(e_{-\alpha}e_{-\beta}\otimes e_{\alpha}e_{\beta}) +
t^\prime_{\alpha\beta}(e_{-\alpha}e_{-\beta}\otimes e_\beta e_\alpha) +
\ldots \bigr)~. \cr} \eqno(2.1)
$$
\no The $H_a$ are generators of the Cartan subalgebra, the $e_\alpha$ are
generators associated with simple roots,
$\varphi^{ab},~t_\alpha,~
t_{\alpha\beta},~t^\prime_{\alpha\beta},~\ldots$ are in the field; the
unwritten terms are monomials in the
$e_\alpha$ and
$e_{-\alpha}$.

More generally, consider the expression (2.1) in the more general case when $H_a,
e_{\pm\alpha}$ generate an associative algebra  with the following relations.
\vskip.5cm

\no {\bf Definition 2.1.}  Let $M,N$ be two countable sets,
$\varphi,\psi$  two maps,
$$
\eqalign{& \varphi :~~M\times M \rightarrow \Crm~, \cr & \psi
:~~M\times N \rightarrow \Crm~, \cr} \quad
\eqalign{a,b &\mapsto \varphi^{ab}~, \cr a,\beta & \mapsto
H_a(\beta)~. \cr} \eqno(2.2)
$$
\vskip.5cm
\no Let ${\cal{A}}$ or ${\cal{A}}(\varphi,\psi)$ be the universal,
associative, unital  algebra over \Crm \enskip with generators
$\{H_a\}\, a\in M,~ \{e_{\pm \alpha}\}\,\alpha \in N$, and relations
$$
\eqalignno{&[H_a,H_b]=0~, \quad  [H_a,e_{\pm\beta}] = \pm
H_a(\beta)e_{\pm\beta}~, & (2.3) \cr
&[e_\alpha,e_{-\beta}]=\delta^\beta_\alpha
\bigl(e^{\varphi(\alpha,\cdot)}-e^{-\varphi(\cdot,\alpha)}\bigr)~, & 
(2.4) \cr}
$$
\no with $\varphi(\alpha,\cdot)=\varphi^{ab}H_a(\alpha)H_b,~
\varphi(\cdot,\alpha)=\varphi^{ab}H_aH_b(\alpha)$ and $
e^{\varphi(\alpha, \cdot) + \varphi(\cdot,\alpha)} \neq 1, ~~ \alpha \in N
$. 

\vskip.5cm
\no The last condition on the parameters is included in order to avoid
having to make some rather trivial exceptions.  

The free subalgebra generated by
$\{e_\alpha\}~\alpha\in N$ (resp.
$\{e_{-\alpha}\}~ \alpha\in N$) will be denoted ${\cal{A}}^+$ (resp.
${\cal{A}}^-$); these subalgebras are  \Zrm ~- graded  , the generators having
grade 1. The subalgebra generated by $\{H_a\}\, a\in M$ is denoted
${\cal{A}}^0$. If necessary we may assume that $M$ is finite.
\vskip.5cm

\no {\bf Definition 2.2.}  A standard R-matrix is a formal series of the
form
$$
\eqalign{R &= \exp\bigl(\varphi^{ab}H_a\otimes H_b\bigr)
\bigl(1+e_{-\alpha}\otimes e_\alpha  + \sum^\infty_{k=2}
t^{\alpha^\prime_1\ldots a^\prime_k}_ {\alpha_1\ldots \alpha_k}
e_{-\alpha_1}\ldots e_{-\alpha_k}\otimes e_{\alpha^\prime_1}\ldots
e_{\alpha^\prime_k}\bigr)~. \cr}
\eqno(2.5)
$$
\vskip.5cm
\no In this formula, and in others to follow, summation over repeated
indices is implied. 
\ve
\no For fixed
$(\alpha)=\alpha_1,\ldots,\alpha_k$ the sum over $(\alpha^\prime)$ runs
over the permutations of $(\alpha)$. The coefficients
$t^{(\alpha^\prime)}_{(\alpha)}$ are in \Crm.
\vskip.5cm

The special property associated with the qualification ``standard" is
that ``$R$ commutes with Cartan"; indeed 
$[R,~H_a\otimes 1+1\otimes H_a]=0,~a\in M$.

We shall determine under what conditions on the parameters 
$\varphi^{ab},~H_a(\beta)$ of ${\cal{A}}$, and for what values of the
coefficients $t^{(a^\prime)}_{(\alpha)}$, the R-matrix (2.5) satisfies
the Yang-Baxter relation
$$ 
Y\hskip-1.0mm B~:=~R_{12}R_{13}R_{23}-R_{23}R_{13}R_{12}=0~. \eqno(2.6)
$$
\no This expression is a formal series in which each term has the form
$\psi_1\otimes \psi_2 \otimes \psi_3\in {\cal{A}}\otimes {\cal{A}}
\otimes {\cal{A}}$.  We assign a double grading as follows.  First extend
the grading of ${\cal{A}}^+$ to the subalgebra of
${\cal{A}}$ that is generated by $\{H_a\} \,a\in M$ and
$\{e_\alpha\}\,\alpha\in N$, by assigning grade zero to $H_a$, and
similarly for
${\cal{A}}^-$. Then $\psi_1$ and $\psi_3$ (but not $\psi_2$) belong to
graded subalgebras of ${\cal{A}}$.  If $\psi_1$ and $\psi_3$ have grades
$\ell$ and $n$, respectively, then define
$$ {\rm grade}~(\psi_1\otimes \psi_2\otimes \psi_3)=(\ell,n)~.
\eqno(2.7)
$$
\no To give a precise meaning to (2.6) we first declare that we mean for
this relation to hold for each grade $(\ell,n)$ separately. This is not
enough, for the number of terms contributing to each  grade is infinite
in general.  The appearance of exponentials in the
$H_a$ can be dealt with in the same way as in the case of simple quantum
groups [TV].  If the sets $M,N$ are infinite, then all results are   
basis dependent.  Eq.(2.6) means that $Y\hskip-1.0mm B$, projected on any
finite subalgebra of ${\cal A}$, vanishes on each grade; the statement
thus involves only finite sums. The analysis of (2.6) will be organized
by ascending grades.
\vskip.5cm
\no {\bf Remarks.}  (i) It is an immediate consequence of (2.6), in grade
(1,1), that
$$ [e_\alpha,e_{-\beta}]=\delta^\beta_\alpha
\bigl(e^{\varphi(\alpha,\cdot)}-e^{-\varphi(\cdot,\alpha)}\bigr)~.
\eqno(2.8)
$$ This relation was therefore included in the definition of the algebra
${\cal A}$.
\no (ii) No relations of the Serre type have been imposed; in fact no 
relations whatever on the subalgebras ${\cal{A}}^+$ and
${\cal{A}}^-$, they are freely generated  by the $e_\alpha$ and by the
$e_{-\alpha}$, respectively. The contextual meaning of  such relations, including
relations of the Serre type, will be discussed in Sections 3-5 and
especially in Section 7.  
\vskip.5cm

 Before stating the main result, it will be convenient to show
the direct evaluation of $Y\hskip-1.0mm B$ up to grade (2,2).  We expand
$$
 R^0~:=~ \exp(\varphi^{ab} H_a\otimes H_b) = R^i\otimes R_i~, \eqno(2.9)
$$
\no sums over $a,b,i$ implied.  Then
$$  e_{-\alpha}R^i\otimes R_i = R^ie_{-\alpha}\otimes 
e^{\varphi(\alpha,\cdot)} R_i~. \eqno(2.10)
$$
\no {\it Grade} (1,1).  The contributions to $R_{12}R_{13}R_{23}$ are of
two kinds:
$$
\eqalign{& R^iR^j e_{-\alpha} \otimes R_iR^k\otimes R_je_\alpha R_k~, \cr
& R^ie_{-\alpha} R^j\otimes R_i e_\alpha R^ke_{-\beta}\otimes
R_jR_ke_\beta~. \cr}
$$
\no Cancellation in $Y\hskip-1.0mm B$ is equivalent to Eq.(2.4).
\vskip 0.5cm
\no {\it Grade} (1,2).  The contributions to $R_{12}R_{13}R_{23}$ are
$$
\eqalign{& R^iR^je_{-\beta}\otimes R_iR^ke_{-\alpha}\otimes R_je_\beta
R_k e_\alpha~, \cr & t^{\alpha'\beta'}_{\alpha\beta}
R^ie_{-\gamma}R^j\otimes R_ie_\gamma R^k e_{-\alpha}e_{-\beta}\otimes
R_jR_ke_{\alpha'} e_{\beta'}~. \cr}
$$
\no Cancellation in $Y\hskip-1.0mm B$ requires that 
$$
\eqalign{t^{\alpha\beta}_{\alpha\beta}&= 
(1-e^{-\varphi(\alpha,\beta)-\varphi(\beta,\alpha)})^{-1}~, \cr
t^{\beta\alpha}_{\alpha\beta} &= -e^{-\varphi(\beta,\alpha)}
t^{\alpha\beta}_{\alpha\beta}~, \quad \alpha\not=\beta~, \cr
t^{\alpha\alpha}_{\alpha\alpha}&= (1+e^{-\varphi(\alpha,\alpha)})^{-1}~.
\cr}
\eqno(2.11)
$$
\no These conditions are necessary and sufficient that the standard
R-matrix satisfy (2.6) up to grade (2,2).

The obstructions to the existence of coefficients $t^{(\alpha^\prime)}_
{(\alpha)}$ such that $Y\hskip-1.0mm B=0$ up to grade (2,2) are therefore
as follows:
$$
\eqalign{& 1+e^{-\varphi(\alpha,\alpha)}=0~~
\hbox{for some $\alpha\in N$}~, \cr &
1-e^{-\varphi(\alpha,\beta)-\varphi(\beta,\alpha)}=0~~ 
\hbox{for some pair $\alpha\not= \beta$}~. \cr} \eqno(2.12)
$$
\no They are typical of obstructions encountered at all grades.

Let
$$ t_{\alpha_1\ldots\alpha_l}:=
t^{\alpha^\prime_1\ldots\alpha^\prime_\ell}_{\alpha_1\ldots\alpha_\ell}\,
e_{\alpha^\prime_1}\ldots e_{\alpha^\prime_\ell}~. \eqno(2.13)
$$

\no {\bf Theorem 2.}  The standard R-matrix (2.5), on
${\cal{A}}$, satisfies the Yang-Baxter relation (2.6) if and only if the
coefficients $t^{(\alpha^\prime)}_{(\alpha)}$ satisfy the following
recursion relation
$$ [t_{\alpha_1\ldots\alpha_\ell},e_{-\gamma}]=
 e^{\varphi(\gamma,\cdot)}
\delta^\gamma_{\alpha_1}t_{\alpha_2\ldots\alpha_\ell}-
t_{\alpha_1\ldots\alpha_{\ell-1}}\delta^\gamma_{\alpha_\ell}e^{-\varphi(\cdot,\gamma)}~.
\eqno(2.14)
$$
\vskip.5cm

\no {\bf Proof.}  (First part.)  We shall prove that (2.14) is 
necessary~--~the ``only if\thinspace" part.  Then we shall study the integrability
of (2.14).  Later, in Section 6, we shall complete the proof of Theorem 2.
Insert (2.5) into (2.6) and use (2.10).  The contribution to
$Y\hskip-1.0mm B$ in grade $(\ell,n)$ is
$$
 R_{12}^0R_{13}^0R_{23}^0e_{-\alpha_1}\ldots e_{-\alpha_\ell} \otimes
P^{\gamma_1\ldots\gamma_n}_{\alpha_1\ldots\alpha_\ell} \otimes
e_{\gamma_1}\ldots e_{\gamma_n}~,
$$
\no in which $P$ is the sum over $m$, from 0 to ${\rm min}(\ell,n)$, of
the following elements of ${\cal A}$,
$$
\eqalign{t_{\alpha_{\ell-m+1}\ldots\alpha_\ell}^{\gamma_1\ldots\gamma_m}
&t_{\alpha_1\ldots\alpha_{\ell-m}}e^{-\varphi(\cdot,\sigma)}
t^{\gamma_{m+1}\ldots\gamma_n} \cr
&-t_{\alpha_1\ldots\alpha_m}^{\gamma_{n-m+1}\ldots\gamma_n}
t^{\gamma_1\ldots\gamma_{n-m}} e^{\varphi(\tau,\cdot)}
t_{\alpha_{m+1}\ldots\alpha_\ell}~, \cr} \eqno(2.15)
$$
\no where $\sigma=\gamma_1+\ldots +\gamma_m$ and
$\tau=\alpha_1+\ldots+\alpha_m$.  The Yang-Baxter relation is satisfied
in grade $(\ell,n)$ if and only if this expression, summed over $m$,
vanishes for every index set $\alpha_1,\ldots,\alpha_\ell$ and
$\gamma_1,\ldots,\gamma_n$.  This is so because ${\cal A}^+$  and ${\cal
A}^-$ are freely generated. We have used the definition (2.13) and
$$ t^{\gamma_1\ldots\gamma_n}~:=~
t_{\gamma^\prime_n\ldots\gamma^\prime_n}^{\gamma_1\ldots\gamma_n}
e_{-\gamma^\prime_1}\ldots e_{-\gamma^\prime_n}~. \eqno(2.16)
$$
\no The lowest grades in which $t_\ell=(t_{\alpha_1\ldots\alpha_\ell})$
appears are $(\ell,0)$ and $(0,\ell)$.  In these cases $m=0$ and (2.15)
vanishes identically.  At grade $(\ell,1)$ one finds (summing
$m=0,1$), the linear recursion relation
$$  [t_\ell,e_{-\gamma}]=e^{\varphi(\gamma,\cdot)}
\delta^\gamma_{\alpha_1}t_{\ell-1}-t_{\ell-1}\delta^\gamma_{\alpha_\ell}
e^{-\varphi(\cdot,\gamma)}~, \eqno(2.17)
$$
\no the full expression for which is Eq. (2.14).  This equation is
therefore necessary.  That it is also sufficient will be proved in
Section 6.
\ve

\line {\bf 3. Differential Algebras. \hfil}
\s

Let $B$ be the unital \Crm-algebra freely generated by
$\{\xi_i\}~i\in  N$, countable.  Suppose given a map
$$ q: N\times N\rightarrow \Crm~, \quad (i,j)\rightarrow q^{ij} \not=0~.
\eqno(3.1)
$$
\no Introduce the natural grading on $B$, $B=\bigoplus B_n$, and a set of
differential operators
$$
\partial_i:~B_n\rightarrow B_{n-1}~,~~i\in N~, \eqno(3.2)
$$
\no defined by
$$
\partial_i\xi_j=\delta_i^j+q^{ij}\xi_j\partial_i~. \eqno(3.3)
$$

We study the problem of integrating sets of equations of the type:
$$
\partial_iX=Y_i~, ~~X\in B~, ~~Y_i\in B~, ~~ i\in N~.   \eqno(3.4)
$$
\no The collection $\{Y_i\}~i\in N$ can be interpreted as the components
of a $B$-valued one-form $Y$, on the space
$\{c^i\partial_i~, ~~ c^i\in \Crm~,~~i\in N\}$.  A constant in
$B_n$ is an element $X\in B_n$, such that $\partial_iX=0, \forall i \in N$.

\vskip.5cm

\no {\bf Proposition 3.1.} (a) The following statements are equivalent: (i)
Eq. (3.4) is integrable for every one-form $Y$ with components in
$B_{n-1}$.  (ii) There are no constants in $B_n$.   (b) When the
parameters $q^{ij}$ are in general position, then there are no constants
in $B_n,~n\geq 1$.
\vskip.5cm

\no {\bf Proof.}  Let
$$ X=X^{i_1\ldots i_n} \xi_{i_1}\dots \xi_n~ \in B_n,
$$
\no then $\partial_iX=0$ means that, for each index set,
$$
\eqalign{X^{i_1\ldots i_n}&+q^{i_1i_2}X^{i_2i_1i_3\ldots}
+q^{i_1i_2}q^{i_1i_3}X^{i_2i_3i_1i_4\ldots} \cr &+\ldots +
q^{i_1i_2}\ldots q^{i_1i_n} X^{i_2\ldots i_ni_1}=0~. \cr} \eqno(3.5)
$$
\no Now fix the unordered index set $\{i_1,\ldots,i_n\}$.  If the values
are distinct then we have a set of $n!$ equations for $n!$ coefficients;
in general the number of equations is always equal to the number of
unknowns.  Solutions exist if and only if the determinant of the matrix
of coefficients vanishes.  This determinant is an algebraic function of
the $q^{ij}$, and not identically zero (Proposition 3.2 below), therefore  
solutions of (3.5), other than $X=0$, exist only on
an algebraic subvariety of parameter space.

The calculation of all these determinants  will be reported elsewhere.
  For $n=2$ the result is
$$D^{12}=1-q^{12}q^{21}~, ~~D^{11}=1+q^{11}~. \eqno(3.6)
$$
\no For $n=3$,
$$
\eqalign{D^{123}&=(1-\sigma^{12})(1-\sigma^{13})(1-\sigma^{23})
(1-\sigma^{12}\sigma^{13}\sigma^{23})~, \cr
D^{112}&=(1+q^{11})(1-\sigma^{12})(1-q^{11}\sigma^{12})~, \cr
D^{111}&=1+q^{11}+(q^{11})^2~, \quad \sigma^{12}~:=~q^{12}q^{21}~. \cr}
\eqno(3.7)
$$
\no It is natural to pass from $B$ to the quotient by the ideal generated
by the constants.  In $B_2$ the constants are
$$
\eqalignno{&\xi_1\xi_2-q^{21}\xi_2\xi_1~, \quad \hbox{when}\quad
\sigma^{12}=1~, & (3.8)\cr &\xi_1\xi_1~~~~~~~~~~~~~~, \quad
\hbox{when}\quad q^{11}=-1~. & (3.9)\cr}
$$
\no If $q^{ii}=-1,~i\in N$ and $\sigma^{ij}=1,~i\not= j$, then the
quotient is a $q$-Grassmann algebra or quantum antiplane.  The constants
in $B_3$ are
$$
\eqalignno{&\xi_1\xi_1\xi_1~~~~~~~~~~~~~~~~~~~~~~~~~~~~~~~~~~~~~~~~~~~~~,
~~1+q^{11}+(q^{11})^2=0~, & (3.10)\cr
&\xi_1\xi_1\xi_2-(q^{21})^2~\xi_2\xi_1\xi_1 ~~~~~~~~~~~~~~~~~~~~~~~~~,
~~1+q^{11}=0~, & (3.11)\cr
&q^{12}\xi_1\xi_1\xi_2-(1+\sigma^{12})~\xi_1\xi_2\xi_1+
q^{21}\xi_2\xi_1\xi_1~,~~ q^{11}\sigma^{12}=1~; & (3.12)\cr}
$$
\no if $\sigma^{12}=1$, there are two constants
$$
\eqalignno{&q^{11}\xi_1\xi_1\xi_2-(1+q^{11})\xi_1\xi_2\xi_1+
(q^{21})^2~\xi_2\xi_1\xi_1~, & (3.13) \cr
&[[\xi_1,\xi_2]_{q^{21}},\xi_3]_{q^{31}q^{32}}~, \quad
[a,b]_q~:=~ab-qba~, & (3.14) \cr}
$$
\no and finally if $\sigma^{12}\sigma^{13}\sigma^{23}=1$ there is one,
$$
\biggl({1\over q^{31}}-q^{13}\biggr)
(\xi_1\xi_2\xi_3+q^{31}q^{32}q^{21}\xi_3\xi_2\xi_1)+\hbox{cyclic.}
\eqno(3.15)
$$
\no Annulment of (3.8), (3.12) and (3.13) are $q$-deformed Serre 
relations [D1].  The last item, Eq. (3.15), may be something new; it
should be interesting to study the quotient of the algebra $B$ with 3
generators by the ideal generated by this constant.

A constant that involves only one variable, $\xi_1$ say, exists if and
only if $q^{11} \neq 1$ is a root of unity,
$$
\xi_1^n ~~{\rm constant ~~ iff} ~~(q^{11})^n = 1, ~~n = 2,3,...,\quad q \neq 1 \,.
$$

It is easy to determine all constants of the $q$-Serre type; that is, all
those that involve two generators and one of them linearly,
$$ C := \sum^k_{m=0} Q^k_m~(\xi_1)^m\xi_2 (\xi_1)^{k-m}~= 0.
\eqno(3.16)
$$
\no With $q=q^{11}$,
$$
\partial_1(\xi_1)^m=(m)_q~(\xi_1)^{m-1}~, \quad (m)_q~:=~1+q+\ldots
+q^{m-1}~. \eqno(3.17)
$$
\no Setting $\partial_1C=0$ gives, for $q^n\not= 1,~n\in Z\hskip-2mm Z$,
$$ Q^k_m=(-q^{12})^m~q^{m(m-1)/2}
\left(\matrix{k\cr m\cr}\right)_q~, \eqno(3.18)
$$
\no while $\partial_2C=0$ is equivalent to
$$
\prod^{k-1}_{m=0}(1-q^m\sigma^{12})=0~.   \eqno(3.19)
$$
\no When $k=2$, compare $D^{112}$ in Eq. (3.7).  If $k$ is the smallest
integer such that a relation like (3.16) holds, then
$$  1-q^{k-1}\sigma^{12}=0~. \eqno(3.20)
$$

Here are some partial results for $B_4$ and $B_5$. $D_{1234}$ is the
product of 12 factors of the form $1 - \sigma_{ij}$, 4 factors of the
form $1-\sigma_{ij}\sigma_{kl}\sigma_{mn}$, 2 identical factors of the
form $1 - \sigma_{12}\ldots\sigma_{34}$; each group accounts for 24
orders in the $q$'s. $D_{12345}$ is the product of 60 factors of the
first type, 20 factors of the second type, 10 factors of  the third type
and 6 identical factors of the form $1 - $(product of all ten
$\sigma_{ij}, i\neq j)$; each group accounts for 120 orders in the $q$'s. Finally the following is true. 
\b
\no {\bf Proposition 3.2.} If all (finite) products of the form $\Pi_{i,j}(q_{ij})^{n_{ij}}$,
where $n_{ij}$ are non-negative integers, differ from unity, then the determinant of 
the matrix of coefficients in (3.5) is different from zero.
 
\ve

\line {\bf 4. Differential Complexes. \hfil}
\s

In the generic case, when there are no constants in $B_n$, the equation
$\partial_iX=Y_i,~Y_i\in B_{n-1},~i\in N$, is always solvable, for any
one-form $Y$.  All one-forms are exact, to be closed has no meaning and
the  differential complex is highly trivial.

The existence of a constant $C\in B_n$ implies that there are one-forms
valued in $B_{n-1}$ that are not exact.  To each 1-dimensional space of
constants in $B_n$ there is a one-dimensional space of non-exact
one-forms, valued in $B_{n-1}$, defined modulo exact one-forms and
obtainable as a limit of $\partial_iX$ as $X\rightarrow C$, after
factoring out a constant.  Thus,
$$ X=\xi_1\xi_2-q^{21}\xi_2\xi_1 \eqno(4.1)
$$
\no becomes a constant as $\sigma^{12}\rightarrow 1$, and a
representative of the associated class of non-exact one-forms is  given by
$$ Y_i=\lim (1-\sigma^{12})^{-1}~\partial_iX =
\cases{\xi_2, & $i=1$, \cr 0~, & $i\not= 1$. \cr} \eqno(4.2)
$$

\vskip 0.5cm
\no {\bf Definition 4.1.} An elementary constant is a linear combination
of re-orderings (permutations of the order of the factors) of a fixed
monomial.
\vskip 0.5cm

A constant $C\in B_n$ also implies a concept of closed one-forms.
\vskip.5cm

\no {\bf Proposition 4.}  If $C\in B_n$,
$$ C=C^{i_1\ldots i_n}\xi_{i_1}\ldots \xi_{i_n}~, \eqno(4.3)
$$
\no is a constant, then the differential operator
$$
\Phi(C)~:=~ C^{i_1\ldots i_n} \partial_{i_1}\ldots\partial_{i_n}
\eqno(4.4)
$$
\no vanishes on $B$.
\vskip.5cm

\no {\bf Proof.}   A constant in $B_n$ is a sum of elementary constants;
it is enough to prove the proposition for the case that $C$ is an
elementary constant.  This implies that there are non-zero $f_i\in \Crm ,
\, i \in N,$ such that the following operator identity 
$$
\partial_iC-f_iC\partial_i=0,\,\, i \in N, \eqno(4.5)
$$
\no holds on $B$.  Let $B^*$ be the unital, associative algebra freely
generated by $\{\partial_i\}~i\in N$, and let
$\Phi:~B\rightarrow B^*$ be the unique isomorphism of algebras such that
$\Phi(\xi_i)=\partial_i$. Let $BB^*(q)$ be the unital, associative
algebra generated by
$\{\xi_i,\partial_i\}~i\in N$, with relations (3.3); then $\Phi$ extends
to a unique isomorphism
$$
\Phi^\prime:~BB^*(q)\rightarrow BB^*(\hat q)~, \quad
\hat q^{ij}=1/q^{ji}~.
$$
\no In particular, $\Phi'(\xi_i) = \partial_i$ and
$\Phi^\prime(\partial_i)=-(q^{ii})^{-1}\xi_i,~i\in N$. Now Eq.(4.5) means
that
$\partial_i\circ C=Cf_i\circ
\partial_i$, where $a\circ b$ denotes the product in $BB^*(q)$.  Applying
$\Phi^\prime$ one gets
$$
\Phi(C) \circ \xi_i=(f_i)^{-1}\xi_i\circ \Phi(C)~,
$$
\no implying that $\Phi(C) X=0,~X\in B$.
\vskip.5cm

\no {\bf Definition 4.2.}  Let $C$ be an elementary constant in
$B_n,~n\geq 2$.  A $B_1^*$ one-form $Y$, valued in $B$, will be said to
be $C$-closed if
$$ d_CY~:=~ C^{i_1\ldots i_n}\partial_{i_1}\ldots\partial_{i_{n-1}}
Y_{i_n}=0~.\eqno(4.6)
$$
\vskip.5cm

\no {\bf Examples.}  In $B_2$ the constants are of the type
$C=\xi_1\xi_1$ or $C^\prime=\xi_1\xi_2-q^{21}\xi_2\xi_1$. Now $Y$ is
$C$-closed if $dY~:=~\partial_1Y_1=0$ and $C^\prime$-closed if
$d^\prime Y~:=~\partial_1Y_2-q^{21}\partial_2Y_1=0$.  The first case is
characteristic of Grassmann algebras and the other of quantum planes. Let
${\cal{C}}$ be the collection
$$ 
\bigl\{C_{ij}=\xi_i\xi_j-q^{ji}\xi_j\xi_i~, \quad i,j\in N~, ~~i\not=
j~\bigr\},\eqno(4.7)
$$
\no and suppose all of them constant.  (In other words, $q^{ij}q^{ji} = 1, i \neq j$.)Then we say that a one-form
$Y$ is ${\cal{C}}$-closed if $Y$ is $C_{ij}$-closed for all
$i\not= j$:
$$
\partial_iY_j-q^{ji}\partial_jY_i=0~, \quad i,j\in N~, ~~i\not= j~.\eqno(4.8)
$$
\no In this case the closure of a $B_1^*$ one-form is expressed in terms
of the $B_1^*$ two-form
$$
Z=dY~, \quad Z_{ij}=\partial_iY_j-q^{ji}\partial_jY_i~,\eqno(4.9)
$$
\no and this naturally leads to familiar $q$-deformed de Rham complexes,
with trivial cohomology. (Non-trivial cohomology depends on completion of the algebra.)

It would be interesting to develop the analogue of this construction, the $q$-deformed de Rham 
complex, in the more general case when ${\cal C}$ is an arbitrary collection of constants. As 
a highly non-trivial example consider the following. Replace the constants $C_{ij}$ in 
(4.7), they are Serre relations of order one, with Serre relations of order two:
$$
C_{ij} = q^{ij}\xi_i\xi_i\xi_j - (1 + \sigma^{ij})\xi_i\xi_j\xi_i + q^{ji}\xi_j\xi_i\xi_i.\eqno(4.10)
$$
This implies that $ q^{ii} = q^{jj} = 1/\sigma^{ij}, \,i \neq j$. Then Definition 4.2 says that 
$Y = (Y_1, Y_2, \ldots)$ is a closed one-form if 
$$
(dY)_{ij} := q^{ij}\partial_i\partial_iY_j - (1 + \sigma^{ij})\partial_i\partial_jY_i
+ q^{ji}\partial_j\partial_iY_i = 0.\eqno(4.11)
$$
Every exact one-form is closed by Proposition 4; the converse statement is less obvious. And 
then there is this problem: what is the integrability condition for the following set of equations
$$
(dY)_{ij} = Z_{ij}, \,\,\, i,j \in N.
$$
In other words, what two-forms are ``closed"?

\vskip1.5cm

\line {\bf 5. Integrability of Eq. (2.14). \hfil}
\s

It was seen, in Section 2, that a necessary condition for the standard
R-matrix (2.5) to satisfy the Yang-Baxter relation (2.6) is that the
coefficients $t^{(\alpha^\prime)}_{(\alpha)}$ satisfy the recursion relation (2.14),
$$
\eqalign{&[t_{\alpha_1\ldots\alpha_\ell},e_{-\gamma}]=
 e^{\varphi(\gamma,\cdot)}\delta^\gamma_{\alpha_1}
t_{\alpha_2\ldots\alpha_\ell}- t_{\alpha_1\ldots\alpha_{\ell-1}}
\delta^\gamma_{\alpha_\ell}e^{-\varphi(\cdot,\gamma)} ~, \cr &
t_{\alpha_1\ldots\alpha_\ell}~:=~
t^{\alpha^\prime_1\ldots\alpha^\prime_\ell}_ {\alpha_1\ldots\alpha_\ell} 
e_{\alpha^\prime_1}\ldots e_{\alpha^\prime_\ell}~. \cr}
\eqno(5.1)
$$
\ve

\no Define \footnote*{ This is where we need the last condition in
Definition 2.1.} two differential operators, $\vec\partial\hskip
-1mm_{-\gamma}$  and $\overleftarrow \partial\hskip -1mm_{-\gamma}$, on
${\cal{A}}^+$, by
$$ [X,e_{-\gamma}]=e^{\varphi(\gamma,\cdot)}\vec\partial_{-\gamma}
X-X\overleftarrow\partial\hskip -1mm_{-\gamma}e^{-\varphi(\cdot,\gamma)}~,
\eqno(5.2)
$$
\no $X\in {\cal{A}}^+$; note that $\overleftarrow\partial\hskip
-1mm_{-\gamma}$ operates from the right.  Similarly,
$$ [e_\alpha,Y]=Y\overleftarrow\partial\hskip -1mm_\alpha
e^{\varphi(\alpha,\cdot)}-e^{-\varphi(\cdot,\alpha)}
\vec\partial_\alpha Y \eqno(5.3) 
$$
\no defines two differential operators on ${\cal{A}}^-$.  These
definitions are equivalent to the rules
$$\matrix{
 \,\vec\partial_{-\gamma} e_\alpha =
\delta^\gamma_\alpha+e^{-\varphi(\gamma,\alpha)} e_\alpha
\vec\partial_{-\gamma}~, &  e_\alpha \overleftarrow\partial\hskip
-1mm_{-\gamma} =
\delta^\gamma_\alpha + e^{-\varphi(\alpha,\gamma)} 
\overleftarrow\partial\hskip -1mm_{-\gamma} e_\alpha, \cr 
e_{-\alpha}\overleftarrow\partial\hskip -1mm_\gamma =
\delta^\gamma_\alpha+e^{-\varphi(\gamma,\alpha)}
\overleftarrow\partial\hskip -1mm_\gamma e_{-\alpha}, &
\,\vec\partial_\gamma e_{-\alpha} =
\delta^\gamma_\alpha+e^{-\varphi(\alpha,\gamma)}
e_{-\alpha}\vec\partial_\gamma~. \cr} \eqno(5.4)
$$
\no Eq. (5.1) is equivalent to
$$
\vec\partial_{-\gamma}t_{\alpha_1\ldots\alpha_\ell}=
\delta^\gamma_{\alpha_1}t_{\alpha_2\ldots\alpha_\ell}~, ~~
t_{\alpha_1\ldots\alpha_\ell}\overleftarrow\partial\hskip -1mm_{-\gamma}=
  t_{\alpha_1\ldots \alpha_{\ell-1}}\delta^\gamma_{\alpha_\ell} ~.
\eqno(5.5)
$$
\vskip.5cm

\no {\bf Proposition 5.1.}  Suppose that the parameters
$\varphi(\alpha,\beta)$ are in general position, so that there are no
constants in
${\cal{A}}^+~({\cal{A}}^-)$ with respect to the differential operators
$\vec\partial_{-\gamma}$ or
$\overleftarrow\partial\hskip -1mm_{-\gamma}$ ($\vec\partial_\gamma$ or
$\overleftarrow\partial\hskip -1mm_\gamma$).  Then either one of the two
equations in (5.5) determines $t_{\alpha_1\ldots\alpha_\ell}$ recursively
and uniquely (the same in each case) from
$t_\alpha=e_\alpha$.
\vskip.5cm

\no {\bf Proof.}  From (5.4) we deduce that
$$ (\vec\partial_{-\gamma}X)
\overleftarrow\partial\hskip -1mm_{-\gamma^\prime}=
\vec\partial_{-\gamma} (X\overleftarrow\partial\hskip
-1mm_{-\gamma^\prime})~. \eqno(5.6)
$$
\no By Proposition 3, the first of (5.5) determines
$t_\ell=t_{\alpha_1\ldots\alpha_\ell}$ uniquely from
$t_\alpha=e_\alpha$.  The other recursion relation also has a unique
solution, $t_\ell^\prime$ say.  We must show that
$t_\ell=t_\ell^\prime,~\ell > 1$.  Since parentheses are superfluous,
$$
\eqalign{\vec\partial_{-\gamma}\,t_\ell\,
\overleftarrow\partial\hskip -1mm_{-\gamma^\prime} &=
\delta^\gamma_{\alpha_1}\,t_{\ell-1}\,
\overleftarrow\partial\hskip -1mm_{-\gamma^\prime}~, \cr
\vec\partial_{-\gamma} \,t^\prime_\ell\,
\overleftarrow\partial\hskip -1mm_{-\gamma^\prime} &=
\vec\partial_{-\gamma} \,t^\prime_{\ell-1}\,
\delta^{\gamma^\prime}_{\alpha_\ell}~. \cr}
$$
\no Suppose $t^\prime_k=t_k$ for $k=1,\ldots,\ell-1$; then the right-hand
sides are both equal to $\delta^\gamma_{\alpha_1}
\,t_{\ell-2}\,\delta^{\gamma^\prime}_{\alpha_\ell}$.  Then the left-hand
sides are also equal and, since there are no constants,
$t_\ell=t^\prime_\ell$.  Since $t_1=t^\prime_1$
($t_\alpha=t^\prime_\alpha=e_\alpha$), the proposition follows by
induction.
\vskip.5cm 
We also encounter the relation
$$
\eqalign{&[e_\alpha,t^{\gamma_1\ldots\gamma_n}]=
t^{\gamma_1\ldots\gamma_{n-1}}\delta_{\alpha }^{\gamma_n}
e^{\varphi(\alpha,\cdot)}-
e^{-\varphi(\cdot,\alpha)}\delta^{\gamma_1}_\alpha
t^{\gamma_2\ldots\gamma_n}~, \cr & t^{\gamma_1\ldots\gamma_n}~:=~
t^{\gamma_1\ldots\gamma_n}_{\gamma_1^\prime\ldots\gamma^\prime_n}
e_{-\gamma^\prime_1}\ldots e_{-\gamma^\prime_n}~. \cr} \eqno(5.7)
$$
\no Just as (5.1) is equivalent to (5.5), this relation is the same as
$$ 
t^{\gamma_1\ldots\gamma_n}\overleftarrow\partial\hskip -1mm_\alpha=
 t^{\gamma_1\ldots\gamma_{n-1}} \delta^{\gamma_n}_\alpha~, ~~
\vec\partial_\alpha t^{\gamma_1\ldots\gamma_n}=
\delta_\alpha^{\gamma_1}t^{\gamma_2\ldots\gamma_n}~. \eqno(5.8)
$$
An adaptation of the proof of Proposition 5.1 shows that either one of these two relations 
determines the same set $t_{\alpha_1,...,\alpha_n}$. Finally, we verify that this unique solution of (5.8)
coincides with the solution of (5.1).
\b
\no {\bf Proposition 5.2.} Fix an unordered set $\{\alpha_1....,\alpha_n\}$ and let $S$ be the matrix 
(with entries in \Crm)~defined by 
the natural pairing between the algebras generated by the $e_\alpha$'s and the $\vec \partial_{-\beta}$'s, 
respectively,
$$
S_{(\alpha)}^{(\beta)} = \vec\partial_{-\beta_n}...\vec\partial_{-\beta_1} \,\,e_{\alpha_1}...e_{\alpha_n},
$$
where $(\alpha), (\beta)$ run over the ordered sets that coincide as unordered set with  
$\{\alpha_1....,\alpha_n\}$. Similarly,
$$
{S'_{(\alpha)}}^{(\beta)} = \vec\partial_{\alpha_n}...\vec\partial_{\alpha_1} \,\,e_{-\beta_1}...e_{-\beta_n}.
$$
Then (a) $S= S'$ and (b) the matrix $t = t_{(\alpha)}^{(\beta)}$ is given by
$$
St = tS' = 1.\eqno(5.9)
$$
\b
\no {\bf Proof.} (a) By inspection: Moving the operator $\vec\partial_{-\beta_1}$ to the right till it encounters
and annihilates $e_{\beta_1}$ produces a factor $\prod_\alpha $exp$ (-\varphi(\beta_1,\alpha))$, where the product runs over 
those $e_{\alpha}$'s that stand to the left of $e_{\beta_1}$. The same factor is is produced by moving $e_{-\beta_1}$
to the left in the expression for $S'$. (b) Iteration of the first of Eq.s (5.5) yields $St = 1$ and 
iteration of the second of (5.8) gives $tS' = 1$. 
  
 \ve

\line {\bf 6. Completion of the Proof of Theorem 2. \hfil}
\s

Suppose that the relations (2.14) are satisfied for $\ell\geq 1$. Now fix
$\ell,n,~\alpha_1,\ldots,\alpha_\ell$ and
$\gamma^1,\ldots,\gamma^n$; we must prove that the expression (2.15),
summed over $m$ from $0$ to ${\rm min}~(\ell,n)$, vanishes.

We begin by calculating the sum over $m=0,1$ (step 1); then we postulate
a formula for the partial sum over $m=0,\ldots,k$ (step $k$). We prove
the formula by induction in $k$, and finally show that the sum vanishes
when $k = {\rm min}~(\ell,n)$.

The term $m=0$ in (2.15) is
$$ [t_\ell,t^n]=
t^{\gamma_1\ldots\gamma_n}_{\gamma_1^\prime\ldots\gamma_n^\prime}
\sum^n_{i=1} e_{-\gamma^\prime_1}\ldots e_{-\gamma^\prime_{i-1}} [t_{\ell
},e_{-\gamma^\prime_i}] e_{-\gamma^\prime_{i+1}}\ldots
e_{-\gamma^\prime_n}~. \eqno(6.1)
$$
\no As in the preceding section we often write $t_\ell,t^n$ for
$t_{\alpha_1\ldots\alpha_\ell}$, $t^{\gamma_1\ldots\gamma_n}$. We shall
gradually make the formulas more schematic so as to bring out their
structure.  By (2.14)
$$
\eqalign{=t^{(\gamma)}_{(\gamma^\prime)} \sum_{i=1}^n
e_{-\gamma^\prime_1}&\ldots \bigl( e^{\varphi(\gamma^\prime_i,\cdot)}
\delta^{\gamma^\prime_i}_{\alpha_1} t_{\ell-1}-t_{\ell-1}
\delta^{\gamma_i^\prime}_{\alpha_\ell} e^{-\varphi(\cdot,\gamma^\prime_i)}
\bigr) \ldots e_{-\gamma^\prime_n}~.\cr} \eqno(6.2)
$$
\no The term $m=1$ is
$$
\eqalign{t^{\gamma_1}_{\alpha_\ell}&t_{\alpha_1\ldots\alpha_{\ell-1}}
e^{-\varphi(\cdot,\gamma_1)}t^{\gamma_2\ldots\gamma_n}
-t^{\gamma_n}_{\alpha_1}t^{\gamma_1\ldots\gamma_{n-1}}
e^{\varphi(\alpha_1,\cdot)}t_{\alpha_2\ldots\alpha_\ell} \cr
&=t_{\ell-1}e^{-\varphi(\cdot,\alpha_\ell)}
\vec \partial_{\alpha_l}t^{\gamma_1\ldots\gamma_n}-
t^{\gamma_1\ldots\gamma_n}\overleftarrow\partial\hskip -1mm_{\alpha_1}
e^{\varphi(\alpha_1,\cdot)}t_{\ell-1} \cr
=t_{\ell-1}t^{(\gamma)}_{(\gamma^\prime)}& \sum
e_{-\gamma^\prime_1}\ldots\delta^{\gamma^\prime_i}_{\alpha_\ell}
e^{-\varphi(\cdot,\gamma^\prime_i)}\ldots e_{-\gamma^\prime_n} 
-t^{(\gamma)}_{(\gamma^\prime)} \sum
e_{-\gamma^\prime_1}\ldots\delta^{\gamma^\prime_i}_{\alpha_1}
e^{\varphi(\gamma^\prime_i,\cdot)}
\ldots e_{-\gamma^\prime_n} t_{\ell-1}~. \cr} 
$$
\no This agrees with (6.2) except for the position of $t_{\ell-1}$, and
the sign.  Thus, adding the contributions $m=0,1$ we obtain
$$
\eqalign{t^{(\gamma)}_{(\gamma^\prime)} \sum_{i<j}
\bigl\{&e_{-\gamma^\prime_1}\ldots\delta^{\gamma_i^\prime}_{\alpha_1}
e^{\varphi(\gamma^\prime_i,\cdot)}e_{-\gamma^\prime_{i+1}}\ldots
[t_{\ell-1},e_{-\gamma^\prime_j}]\ldots e_{-\gamma^\prime_n} \cr
&+e_{-\gamma^\prime_1}\ldots[t_{\ell-1},e_{-\gamma^\prime_i}]
e_{-\gamma^\prime_{i+1}}\ldots\delta^{\gamma'_j}_{\alpha_l}
e^{-\varphi(\cdot,\gamma^\prime_j)}\ldots e_{-\gamma^\prime_n}\bigr\}~.
\cr}
\eqno(6.3)
$$
\no This completes step 1; all terms involving $t_\ell$ have disappeared
and $t_{\ell-1}$ appears only in commutators that allow us to use (2.14)
again.

We claim that after carrying out step $k$, which includes summing over
$m=0,\ldots,k$, one obtains the following expression
$$
\sum^k_{s=0} \ldots(\delta e^\varphi)^{k-s}\ldots[t_{\ell-k},e_{-\gamma}]
\ldots (\delta e^{-\varphi})^s \ldots~,\,\, k < {\rm min} (l,n),
\eqno(6.4)
$$ and zero, $k = \,\,$min$ (l,n)$. Here the dots stand for products of
the $e_{-\gamma^\prime_i}$, interrupted $k-s$ times by a factor of the
type
$\delta^{\gamma^\prime_i}_ {\alpha_1}
e^{\varphi(\gamma_i^\prime,\cdot)}$, once by $[~,~]$ and $s$ times by a
factor like
$\delta_{\alpha_l}^{\gamma_i}  e^{-\varphi(\cdot,\gamma^\prime_i)}$, as
in (6.3).

To verify this claim we carry out the next step.  We first evaluate the
commutators and examine the cancellations that take place between
successive terms in the sum (6.4):
$$
\eqalign{\ldots&[t_{\ell-k},e_{-\gamma^\prime_i}] \ldots  (\delta
e^{-\varphi}) \ldots  + \ldots (\delta e^\varphi) \ldots
[t_{\ell-k},e_{-\gamma^\prime_i}] \dots \cr &\quad =
\ldots(\delta^{\gamma^\prime_i}_\alpha e^{\varphi(\alpha,\cdot)}
t_{\ell-k-1}-t_{\ell-k-1}\delta^{\gamma^\prime_i}_{\alpha^\prime}
e^{-\varphi(\cdot,\alpha^\prime)}) \ldots(\delta e^{-\varphi})
\ldots \cr &\quad \quad \quad  +(\delta e^\varphi)\ldots
(\delta^{\gamma^\prime_j}_\alpha
e^{\varphi(\alpha,\cdot)}t_{\ell-k-1}-t_{\ell-k-1}
\delta^{\gamma^\prime_j}_{\alpha^\prime}
e^{-\varphi(\cdot,\alpha^\prime)}) \ldots~. \cr}      
$$
\no The first term in the first line combines with the second term in the
second line to
$$
\ldots (\delta e^\varphi)\ldots[t_{\ell-k-1},e_{-\gamma}] \ldots (\delta
e^{-\varphi}) \ldots~.
$$
\no Successive terms in (6.4) all combine in this way, to reproduce the
same expression with $k$ replaced by $k+1$, except for the fact that
there is no term in the sequence that precedes and collaborates with the
first term and no term that succeeds and collaborates with the last
term.  It remains, therefore, to be proved that the summand
$m=k+1$ in (2.15) precisely supplies the two missing terms.

By (5.8),
$$
\eqalign{& \vec \partial_\beta t^{\gamma_1\ldots\gamma_n} =
\delta^{\gamma_1}_\beta t^{\gamma_2\ldots\gamma_n}~, \cr &
\vec\partial_{\beta_m}\ldots \vec\partial_{\beta_1} t^{\gamma_1\ldots
\gamma_n}=\delta^{\gamma_1}_{\beta_1}\ldots
\delta^{\gamma_m}_{\beta_m} \,t^{\gamma_{m+1}\ldots\gamma_n}~, \cr &
t^{\beta_1\ldots\beta_m}_{\alpha_{\ell-m+1}\ldots\alpha_\ell}
\vec\partial_{\beta_m}\ldots\vec\partial_{\beta_1}
t^{\gamma_1\ldots\gamma_n} =
t^{\gamma_1\ldots\gamma_m}_{\alpha_{\ell-m+1}\ldots\alpha_l}
t^{\gamma_{m+1}\ldots\gamma_n}~. \cr}
$$
\no Hence, if $\vec t_{\alpha_1\ldots\alpha_\ell}$ is the differential
operator
$$
\vec t_{\alpha_1\ldots\alpha_\ell}~:=~
t^{\alpha_1^\prime\ldots\alpha^\prime_\ell}_{\alpha_1\ldots\alpha_\ell}
\vec\partial_{\alpha^\prime_l} \ldots \vec\partial_{\alpha^\prime_1}~,
$$
\no then
$$ t^{\gamma_1\ldots\gamma_m}_{\alpha_{\ell-m+1}\ldots\alpha_\ell}
t^{\gamma_{m+1}\ldots\gamma_n}=
\vec t_{\alpha_\ell\ldots\alpha_{\ell-m+1}} t^{\gamma_1\ldots\gamma_n}
\eqno(6.5)
$$
\no and the first of the two  terms in (2.15) is
$$
\eqalign{t_{\alpha_1\ldots\alpha_{\ell-m}} &e^{-\varphi(\cdot,\sigma)}
\vec t_{\alpha_\ell\ldots\alpha_{\ell-m+1}} t^{\gamma_1\ldots\gamma_n} \cr
&=t_{\alpha_1\ldots\alpha_{\ell-m}}t^{(\gamma)}_{(\gamma^\prime)}
\sum^n_{i=1} e_{-\gamma^\prime_1}\ldots
\bigl[e^{-\varphi(\cdot,\sigma)}
\vec t_{\alpha_\ell\ldots\alpha_{\ell-m+1}}, e_{-\gamma^\prime_i}\bigr]
\ldots \cr  &= t_{\alpha_1\ldots\alpha_{\ell-m}}
t^{(\gamma)}_{(\gamma^\prime)}
\sum_i e_{-\gamma^\prime_1} \ldots 
\delta^{\gamma^\prime_i}_{\alpha_{\ell-m+1}}
e^{-\varphi(\cdot,\sigma)}\vec
t_{\alpha_\ell\ldots\alpha_{\ell-m+2}}\ldots~. \cr}
$$
\no By iteration of these steps one finally ends up, when $m=k+1$, with
precisely the missing terms; actually one of the missing terms, we leave
it to the reader to carry out the calculation for the other one.  This
done, the proof of Theorem 2 is complete.
\vskip.3cm

{\bf Corollary 6.} With the parameters in general position, there exists
a unique standard R-matrix on ${\cal A}$ that satisfies the Yang-Baxter relation.

\vskip1.5cm

\line {7. \bf Obstructions and Generalized Serre Relations. \hfil}

We have been concerned with the construction of a standard R-matrix,
Definition (2.2), that satisfies the Yang-Baxter relation,  Eq.(2.6),  
on an algebra
${\cal{A}}$, Definition (2.1). The relations of ${\cal{A}}$ involve
parameters; when these parameters are in general position then the recursion relation
(2.14) has a unique solution that provides the unique standard R-matrix on ${\cal A}$
that satisfies the Yang-Baxter relation. At certain hypersurfaces in parameter space we have
encountered obstructions, characterized by the vanishing of one or more
of the determinants that we have studied in Section 3.  At these points
there appear elements in ${\cal{A}}^+$ that are constants with respect to
differential operators
$\vec\partial_{-\alpha}$ and $\overleftarrow\partial\hskip -1mm
_{-\alpha}$, and elements in ${\cal{A}}^-$ that are constants with
respect to
$\vec\partial_\alpha$ and $\overleftarrow\partial\hskip -1mm_\alpha$.
Then there is no solution of (2.14) and no standard R-matrix on ${\cal A}$ that 
satisfies Yang-Baxter.

We shall show that all these obstructions can be overcome by the
introduction of additional relations in the definition of ${\cal{A}}$ or,
what is the same, by replacing ${\cal{A}}$ by a quotient ${\cal{A}}/I$,
where $I$ is the ideal generated by the constants.  The next three
propositions relate the null-spaces of the four differential operators to
each other.
\vskip.5cm

\no {\bf Proposition 7.1.}  The space of constants with respect to
$\vec\partial_{-\gamma}$ in ${\cal{A}}^+_n$ has the same dimension as the
space of constants with respect to 
$\overleftarrow\partial\hskip -1mm_{-\gamma}$.  If there are no constants
in
${\cal{A}}^+_\ell$ for $\ell=1,\ldots,n-1$, then the two spaces coincide.
\vskip.5cm

\no {\bf Proof.}  An easy consequence of Eq. (5.6).
\vskip.3cm
\ve
Let $C\in {\cal{A}}^+_n$ be a constant with respect to
$\vec\partial_{-\gamma}$, 
$\gamma \in N$, and assume, provisionally, that there are no constants in
${\cal{A}}^+_\ell,~1<\ell<n$. Without essential loss of generality we
take $C$ to be a linear combination
$$ C=C^{\gamma_1\ldots\gamma_n} e_{\gamma_1}\ldots e_{\gamma_n}~,
$$
 
\no where the summation runs over the permutations of a fixed set
$\{\gamma_1\ldots\gamma_n\}$, hence over a finite set.  Let $d$ be the
operator that takes $X\in {\cal{A}}^+_n$ to the one-form $Y$ valued in
${\cal{A}}^+_{n-1}$ with components $\vec\partial_{-\gamma}X,~\gamma\in
N$. This operator is represented by a direct sum of finite square
matrices, also denoted $d$.  The constant $C$ is a null-vector for $d$. 
The transposed matrix   also has a null-vector; it exists by virtue of
the fact that $dX$ is
$C$-closed: (Definition 4.2):
$$
\eqalign{(C^{\gamma_1\ldots\gamma_n}&
\vec\partial\hskip -1mm_{-\gamma_1}\ldots\vec\partial\hskip
-1mm_{-\gamma_n}) e_{\alpha_1}\ldots e_{\alpha_n} \cr &=
(C^{\gamma_1\ldots\gamma_n}\vec\partial\hskip -1mm_{-\gamma_1}
\ldots\vec\partial\hskip -1mm_{-\gamma_{n-1}})
(d^{\gamma_n\beta_1\ldots\beta_{n-1}}_{\alpha_1\ldots\alpha_n}
e_{\beta_1}\ldots e_{\beta_{n-1}})
\cr &= d^{\gamma_n\beta_1\ldots\beta_{n-1}}_{\alpha_1\ldots\alpha_n}
(C^{\gamma_1\ldots\gamma_n}\vec\partial\hskip -1mm_{-\gamma_1}\ldots
\vec\partial\hskip -1mm_{-\gamma_{n-1}} e_{\beta_1}\ldots
e_{\beta_{n-1}}) =0~. \cr}
$$
\no The obstruction to solving Eq. (5.5) is that the right-hand side is
not in the null-space of the transpose of $d$; it is not
$C$-closed. Indeed, since there are no constants in
${\cal{A}}^+_\ell,~\ell<n$,
$$ (C^{\gamma_1\ldots\gamma_n}\vec\partial\hskip -1mm_{-\gamma_1}\ldots
\vec\partial\hskip -1mm_{-\gamma_{n-1}}) \delta^{\alpha_1}_{\gamma_n}
t_{\alpha_2\ldots\alpha_n}=C^{\alpha_n\ldots\alpha_1}\not= 0~.
$$
\no Recall that the R-matrix is expressed in terms of
$e_{-\alpha_1}\ldots e_{-\alpha_n}t_{\alpha_1\ldots\alpha_n}$. The
obstruction to Yang-Baxter is thus
$$ e_{-\alpha_1}\ldots e_{-\alpha_n}C^{\alpha_n\ldots\alpha_1}~=:~
C^\prime \in {\cal{A}}^-_n~.
$$
\vskip.5cm

\no {\bf Proposition 7.2.}  The element $C^\prime \in {\cal{A}}^-_n$ is a
constant.
\vskip.5cm

\no {\bf Proof.}  One verifies directly that 
$\vec\partial_{-\gamma} C=0,~\gamma \in N$, is equivalent to
$C^\prime\overleftarrow\partial\hskip -1mm_\gamma =0,~\gamma\in N$.
\vskip.5cm

Thus, if the first obstruction to the Yang-Baxter relation is 
encountered at the evaluation of $t_{\alpha_1\ldots\alpha_n}$, then this
obstruction can be avoided by   replacing ${\cal{A}}$ by the quotient
${\cal{A}}/I_n$, where $I_n$ is the ideal generated by the constants in
${\cal{A}}^\pm_n$.  Once this has been done, we study the obstructions at
the next level.  Since the constants at level
$n$ have been removed there are none in ${\cal{A}}^\pm_\ell,~
\ell\leq n$, and substantially the same analysis applies to constants in
${\cal{A}}^\pm_{n+1}$.  To formulate the final result we need:
\vskip.5cm

\no {\bf Proposition 7.3.}  The ideal $I^+ \subset {\cal{A}}^+$ generated
by the constants of $\vec\partial_{-\gamma}$, coincides with that
generated by the constants of
$\overleftarrow\partial\hskip -1mm_{-\gamma}$.  The same statement holds
true in ${\cal{A}}^-$, {\it mutatis mutandi}.
\vskip.5cm

When the parameters are in general position there are no constants, and
Theorem 2 with Proposition 5 assures us that there is a unique standard
R-matrix in ${\cal{A}}\otimes {\cal{A}}$ that satisfies the Yang-Baxter
relation  (2.6).  We are now in a position to allow for the appearance of
constants.

\vskip.5cm

\no {\bf Remark.}  There are no constants in ${\cal{A}}_1^\pm$; the
generators $H_a,e_{\pm\alpha}$ of ${\cal{A}}$ are also generators of
${\cal{A}}^\prime={\cal{A}}/I$.
\vskip.5cm

\no {\bf Theorem 7.}  Let $I\subset {\cal{A}}$ be the ideal generated by
the constants in ${\cal{A}}^+$ and the constants in
${\cal{A}}^-$, and let ${\cal{A}}^\prime$ be the quotient ${\cal{A}}/I$. 
Interpret the standard, universal R-matrix (2.5) as an element of
${\cal{A}}^\prime\otimes {\cal{A}}^\prime$.   The Yang-Baxter relation
for the standard R-matrix on $\cal A'$ is equivalent to the recursion
relation
$$
\eqalign{[t_l,1 \otimes e_{-\gamma}] &= (e_{-\gamma} \otimes
e^{\varphi(\gamma,\cdot)})t_{l-1} - t_{l-1}(e_{-\gamma} \otimes 
e^{-\varphi(\cdot,\gamma)}),\cr t_l &:=
t^{(\alpha')}_{(\alpha)}e_{-\alpha_1} \ldots e_{-\alpha_l} \otimes
e_{\alpha_1'}\ldots e_{\alpha_l'}, \cr } \eqno(7.1)  
$$ and to either one of the following
$$ [e_\gamma,t_l \otimes 1] = t_{l-1}(e^{\varphi(\gamma,\cdot)} \otimes e_\gamma) -
(e^{-\varphi(\cdot,\gamma)} \otimes e_\gamma)t_{l-1},
\eqno(7.2)
$$
$$ (1 \otimes \vec\partial_{-\gamma})t_l= (e_{-\gamma} \otimes 1)t_{l-1},
\quad t_l(1 \otimes \overleftarrow\partial\hskip -1mm _{-\gamma}) =
t_{l-1}( e_{-\gamma} \otimes 1),\eqno(7.3)
$$
$$ t_l(\overleftarrow\partial\hskip-1mm _\gamma \otimes 1) = t_{l-1}(1
\otimes e_\gamma),\quad (\vec\partial_\gamma \otimes 1)t_l =  (1  \otimes
e_\gamma)t_{l-1}.\eqno(7.4)
$$ These relations are integrable (with $t_1 = e_{-\alpha} \otimes
e_{\alpha}$) and yield a unique standard R-matrix on $\cal A'$.
\vskip1.5cm

\ve

\line {\bf 8. The Standard Classical r-matrix for Simple Lie Algebras. \hfil}
\s

We shall now specialize, by stages, until we arrive at simple quantum
groups, where a limiting process relates the standard R-matrix to a classical
r-matrix.

Suppose that 
$$ 
{\rm Card}(N)~:=~\ell<\infty.\eqno(8.1)
$$ 
Suppose next that the
ideal $I$ (generated by the constants of
${\cal{A}}$) is generated by a complete set of Serre relations; that is,
for each pair
$\alpha,\beta\in N,~\alpha\not= \beta$, there is a smallest positive
integer $k_{\alpha\beta}$ such that there is a relation in 
${\cal{A}}/I$ of the form
$$
\sum^{k_{\alpha\beta}}_{m=0} Q_m^{(\alpha,\beta)}(e_\alpha)^m
e_\beta(e_\alpha)^{k_{\alpha\beta}-m} =0~, \eqno(8.2)
$$
\no with coefficients $Q^{(\alpha\beta)}_m$ in the field.  The left side,
as an element of ${\cal{A}}^+$, is a constant, and the  penultimate
paragraph of Section 3 applies.  In particular, the relation (3.20)
becomes
$$ 
e^{\varphi(\alpha,\beta)+\varphi(\beta,\alpha)+(k_{\alpha\beta}-1)\varphi(\alpha,\alpha)}
=1~, \eqno(8.3)
$$
and the coefficients are   
 $$
 Q^k_m=(-)^m e^{m\varphi(\alpha,\beta)}~q^{m(m-1)/2}
\left(\matrix{k\cr m\cr}\right)_q~,\quad q := e^{\varphi(\alpha,\alpha)}.\eqno(8.4)  
$$ 

\no We specialize further by supposing that the exponent in (8.3) vanish,
$$ 
\varphi(\alpha,\beta)+\varphi(\beta,\alpha)=(1-k_{\alpha\beta})
\varphi(\alpha,\alpha)~, \quad \alpha\not= \beta~. \eqno(8.5)
$$ The form $(\cdot,\cdot)$ defined by
$$ (\alpha,\beta)=\varphi(\alpha,\beta)+\varphi(\beta,\alpha) \eqno(8.6)
$$
\no will be called the restricted Killing form, and the $\ell$-by-$\ell$
matrix with components
$$ A_{\alpha\beta}~ = {2(\alpha,\beta) \over
(\alpha,\alpha)} \eqno(8.7)
$$
\no will be called the generalized Cartan matrix; note that it is
symmetrizable. Finally, a suitable restriction on Card($M$) brings us
to   quantized Kac-Moody algebras.

Let ${\cal A}'_{cl}$ be the algebra obtained from
${\cal A}'$ when the relations (2.4) are replaced by
$$ 
[e_\alpha,e_{-\beta}] = \delta_\alpha^\beta\,\bigl( \varphi(\alpha,\cdot) +
\varphi(\cdot,\alpha)\bigr).\eqno(8.8)
$$
 If ${\cal A}_{cl}'$ is a Kac-Moody algebra of finite type, resp. affine type, then we 
may say that ${\cal A}'$ is a quantized Kac-Moody algebra of finite type, resp. affine type.
But because ${\cal A}'$ cannot be recovered from ${\cal A}_{cl}'$ an autonomous definition
is preferable.
\vskip.50cm
\no {\bf Definition 8.} Let ${\cal A'}$ be as above; that is, the quotient of an algebra
${\cal A}$ as per Definition 2.1, with parameters satisfying (8.3), by the ideal  generated
by the Serre relations (8.2). We shall say that ${\cal A}'$  
 is a quantized Kac-Moody algebra of finite type if (i) Card $M$ =   Card $N = l
<\infty$,  and (ii) the (symmetrizable) generalized Cartan matrix
$$
A _{\alpha\beta} =  {\varphi(\alpha,\beta) + \varphi(\beta,\alpha) 
\over \varphi(\alpha,\alpha)}\eqno(8.9)
$$ 
is positive definite with $A_{\alpha\beta} \in \{0,-1,\ldots\}, \,\, \alpha \neq \beta$.
  We shall say that
${\cal A}'$  
 is a quantized Kac-Moody algebra of affine type if (i) Card $M$ = 1 +  Card $N   <
\infty$, and (ii) the  generalized Cartan matrix
 is positive semi-definite with $A_{\alpha\beta} \in \{0,-1,\ldots\}, \,\, \alpha \neq
\beta$ and
  all its principal minors are  positive definite.
\vskip.5cm
\no The remainder of this section deals with Kac-Moody algebras of finite type.

The semi-classical limit of $R$ is defined by replacing
$$
\eqalign{&\varphi(\cdot,\cdot)\rightarrow \hbar\varphi(\cdot,\cdot)~,  
\cr & e_\alpha \rightarrow \kappa \,e_\alpha~, \quad e_{-\alpha}\rightarrow \kappa'
e_{-\alpha}~,
\quad \kappa\kappa^\prime = \hbar~, \quad \alpha\in N~, \cr} \eqno(8.10)
$$
\no and developing the exponentials to first order in $\hbar$.  Then Eq.
(2.4) becomes
$$
\eqalign{[e_\alpha,e_{-\beta}]&= 
\delta_{\alpha\beta}\bigl(\varphi(\alpha,\cdot)
+\varphi(\cdot,\alpha)\bigr)  \cr &=:~\delta_{\alpha\beta}H_{(\alpha)}
\varphi(\alpha,\alpha)~. \cr}
\eqno(8.11)
$$
 (Definition of $H_{(\alpha)} \in \Cit$.) It follows from (8.11) and (2.3) that
$$ [H_{(\alpha)},e_\beta]=A_{\alpha\beta}e_\beta~, \quad
\alpha,\beta\in N~. \eqno(8.12)
$$ The definition (8.7) of the generalized Cartan matrix implies that
$A_{\alpha\alpha}=2,~\alpha\in N$, that $A_{\alpha\beta}\in
\{0,-1,-2,\ldots\},~\alpha\not= \beta$, and that $A_{\alpha\beta}\not= 0$
implies $A_{\beta\alpha}\not= 0$.  Special cases are affine Lie algebras
and simple Lie algebras.  The latter are characterized by two additional
properties of $(A_{\alpha\beta})$: indecomposability and
${\rm det}(A)>0$.  We  now assume that both hold, and that
$\{H_{(\alpha)},\, \alpha \in N\}$ generates ${\cal A}^0$.

The (classical) r-matrix $r$ associated with the standard R-matrix (2.5)
is defined by
$$  R=1+\hbar r+o(\hbar^2)~. \eqno(8.13)
$$
\no Two terms in r are obvious: $r=\varphi  + \sum e_{-\alpha}
\otimes e_\alpha+{\rm ?}$, with the sum extending over simple roots.
Evaluating the remaining terms is more difficult, because a) we do not have a sufficiently explicit  expression for 
the coefficients $t_{(\alpha)}^{(\alpha')}$ and b) because all these coefficients are singular 
in the classical limit. Both these difficulties are avoided by the recursion relation (7.1), 
as we shall see later. The result, which was known by indirect means,
   with a particular normalization of the non-simple roots, is that
$$
 r=\varphi+\sum_{\alpha\in\Delta^+} E_{-\alpha}\otimes E_\alpha~,
\eqno(8.14)
$$
\no where $\Delta^+$ is the set of positive roots.  (Definition 9.) This is what we call the 
the standard r-matrix for a simple Lie algebra.  It satisfies the classical Yang-Baxter relation
$$
 [r_{12},r_{13} + r_{23}]+[r_{13},r_{23}]=0 \eqno(8.15)
$$
\no and
$$
 r+r^t=\hat K~ ,\eqno(8.16)
$$ the Killing form of \g. In the list of (constant) r-matrices
obtained by Belavin and  Drinfeld [BD], (8.14) is the simplest.  The
quantum groups to which these r-matrices are associated are the twisted
quantum groups of Reshetikhin and others [R][Sc][Su].

\vskip1.5cm

 \ve

\no {\bf 9. The Standard Classical r-matrix for Untwisted Loop Algebras. }

A quantized affine Kac-Moody algebra can be described as follows. 
Let $\hat {\cal A}'$ be as above, with parameters satisfying (1.8) and Serre relations 
(1.9), with root generators $\{e_{\pm\alpha}\}\,\alpha = 0,\ldots,l $ and Cartan
generators 
$H_1,\ldots ,H_l,\,c,\,d$, such that the subset that consists of $\{e_{\pm \alpha}\}\,
\alpha \neq 0$ and $H_1, \ldots , H_l$  
 generates a subalgebra ${\cal A}'$ that is a quantized Kac-Moody algebra of
finite type. Let $\hat \varphi$ refer to $\hat {\cal A}'$ and $\varphi$ to ${ \cal A}'$, 
and suppose that
$$
\hat \varphi = \varphi + u \, c\otimes d + (1-u)\, d \otimes c,\quad 
[d, e_{\pm \alpha}] = \pm \delta_\alpha^0 \,e_{\pm 0}.
$$
with some $ u \in \Cit$. 
 Suppose that $c$ is central and that the extra root defined by
$
[H_a,e_0] = H_a(0)e_0
$
is such as to make the generalized Cartan matrix of $\hat {\cal A}'$  positive
semi-definite  with all its principal minors positive. Then $\hat {\cal A}'$ is a 
quantized affine Kac-Moody algebra.

Consider a quantized affine Kac-Moody algebra $\hat{\cal A}'$, with generators $e_{\pm 0},
\ldots ,e_{\pm l}$ and $H_1,\ldots ,H_l, c, d$. Renormalize as in (8.10) and pass to the
classical limit.  
\vskip.5cm
\no {\bf Definition 9.} Positive root vectors are elements in ${\cal A}_{cl}'^+$ defined 
recursively.    (a) The generators $e_\alpha$ are positive root vectors. (b) If 
$E_i$ and $ E_j$ are positive root vectors and $[E_i,E_j] \neq 0$, then  $[E_i,E_j]$ is a positive root vector. 
 (c) All positive root vectors  are obtained in this way from the generators. Negative
root vectors are in ${\cal A }_{cl}'^-$ and are defined analogously.
\vskip.5cm
 
Let $\{E_i\}\,i = 1,\ldots ,n,+$ be
the  positive root vectors, labelled in such a way that
 $$ [e_\alpha,E_+] = 0 = [e_{-\alpha},E_-],\eqno(9.1)
$$ 
and
\vskip-4mm 
$$ 
[E_i, E_-]    \in {\cal A }_{cl}'^0 \cdot {\cal A }_{cl}'^-, 
\quad [E_{-i}, E_+]   \in {\cal A}_{cl}'^0 \cdot {\cal A }_{cl}'^+.\eqno(9.2)
$$ 
Then we may refer to $E_+$ as a highest root vector.

Suppose that the extra root $H_a(0) = H_a(E_-)$, and  
 pass to the associated untwisted loop algebra ~$ \Cit 
[\lambda,\lambda^{-1}] \otimes
{\cal A}_{cl}'$ by substituting
$$
\hat\varphi \rightarrow \varphi,\,\,\,e_0 := \lambda E_-,\quad e_{-0} := \lambda^{-1}E_+.
\eqno(9.3) 
$$
(Replacing $\hat \varphi$ by $ \varphi$ amounts to taking the quotient by the ideal generated by
the central element $c$.)

After the renormalization (8.10) $t_n$ is of order $\hbar$ and the classical r-matrix 
is defined by (8.13),
 $$ 
R = 1 + \hbar r + o(\hbar^2).\eqno(9.4)
$$ 
 
The Yang-Baxter relation for $R$ is equivalent to the  recursion relation (7.1),
$$  
[t_n, 1 \otimes e_{-\gamma}] = (e_{-\gamma} \otimes e^{\varphi(\gamma,\cdot)}) t_{n-1}
 - t_{n-1}   (e_{-\gamma} \otimes e^{-\varphi(\cdot,\gamma)}),\,\, n \geq 1.  \eqno(9.5)
$$ 
To lowest order in $\hbar$ this becomes
$$
\eqalign{
[t_1, 1 \otimes e_{-\gamma}] &= e_{-\gamma} \otimes (\varphi + \varphi^t)(\gamma),\cr 
[t_n, 1 \otimes e_{-\gamma}] &= [e_{-\gamma} \otimes 1,t_{n-1}],\,\,\, 
  n \geq 2,\cr}\eqno(9.6)
$$ 
which is the same as
$$ 
[1 \otimes e_{-\gamma} + e_{-\gamma} \otimes 1,r - \varphi] +  [t_1,1 \otimes e_{-\gamma}  ]
= 0,\,\, \gamma = 0,\cdots,l,\eqno(9.7)
$$ 
with $t_1 = \sum e_{-\alpha} \otimes e_\alpha$, or
$$ 
[1 \otimes e_{-\gamma} + e_{-\gamma} \otimes 1,\,\,r] = \varphi(\cdot,\gamma) \wedge 
e_{-\gamma}.
$$ 
This result is just the classical limit of the relation 
$\Delta(e_{-\gamma})R = R\Delta'(e_{-\gamma})$, which explains why it determines $r$.

We normalize the root vectors so that the Casimir element takes the form
$$
C = \varphi + \varphi^t + \sum E_{-i} \otimes E_i + \sum E_i \otimes E_{-i}.\eqno(9.8)
$$ Then
$$ [e_{-\gamma}, E_{-i}] = cE_{-j} {\rm ~ implies ~ that ~ } [E_{j},e_{-\gamma}] = cE_i,
\,\, \gamma \neq 0, \eqno(9.9)
$$ 
$$ [e_{-0}, E_{-i}] = cE_{j} {\rm ~ implies ~ that ~ } [E_{-j},e_{-0}] = cE_i, \,\,\gamma
\neq 0,\eqno(9.10)
$$  
It may be seen from the structure of $t_n$ that it is a polynomial of order $n$
in $\lambda/\mu$. The recursion relation shows that the classical limit is in 
${\cal A}_{cl}' \otimes {\cal A}_{cl}'$. The classical r-matrix can therefore be expressed as a formal power 
series in $x = \lambda/\mu$,
$$ r = \varphi + \psi(x)^{ab}H_a \otimes H_b + \sum f_i(x) E_{-i} \otimes E_i + \sum g_i(x) E_i
\otimes E_{-i}.\eqno(9.11)
$$  Now it is easy to work out the implications of Eq.(9.7), namely, first taking $\gamma
\neq 0$,
$$
\eqalign{
 0 &=  [1 \otimes e_{-\gamma} + e_{-\gamma} \otimes 1, \psi(x)^{ab}H_a \otimes H_b + \sum
f_i(x) E_{-i} \otimes E_i \cr &\hskip1in + \sum g_i(x) E_i \otimes E_{-i}] + \sum
e_{-\alpha}\otimes [e_\alpha,e_{-\gamma}]\cr 
&= e_{-\gamma} \otimes \bigl(
\psi(\gamma,\cdot) + (1-f_\gamma) (\varphi + \varphi^t)(\gamma) \bigr) \cr 
&\quad \quad +
\bigl(\psi(\cdot,\gamma) - g_\gamma (\varphi + \varphi^t)(\gamma)\bigr) \otimes e_{-\gamma}\cr
&\quad \quad+
\sum f_i[e_{-\gamma},E_{-i}] \otimes E_i + {\sum}' f_iE_{-i} \otimes [e_{-\gamma},E_i]\cr
 &\quad \quad + {\sum}'
g_i[e_{-\gamma},E_i] \otimes E_{-i} + \sum g_i E_i \otimes [e_{-\gamma},E_{-i}], \,\,\gamma
\neq 0.
\cr}\eqno(9.12)
$$ 
The prime on ${\sum}'$ means that the summation is
over  roots that are not simple. Cancellation  in the last two lines  imply, in view of
(9.9) and since the adjoint action is irreducible, that $f_i = f,\,g_i = g, \, i =
1,\ldots,l$. Cancellation in the two first lines now tells us that $\psi \propto 
\varphi + \varphi^t$, hence $\psi$ is symmetric, and it follows that $ g = f-1$.  This
gives us 
$$
r = \varphi + \sum E_{-i} \otimes E_i + g(x) C, \eqno(9.13)
$$
which is actually obvious: The two first terms is a special solution and the last term is
the only thing that commutes with $\Delta_0(e_{-\gamma}) = 1 \otimes e_{-\gamma} +
e_{-\gamma} \otimes 1$. Next, Eq.(9.7) with $\gamma = 0$,
$$
\eqalign{
 0 &=  [1 \otimes e_{-0} + e_{-0} \otimes 1,\, \psi(x)^{ab}H_a \otimes H_b + \sum
f_i(x) E_{-i} \otimes E_i \cr 
&\hskip1in + \sum g_i(x) E_i \otimes E_{-i}] + \sum
e_{-\alpha}\otimes [e_\alpha,e_{-0}]\cr 
&= E_+ \otimes \bigl({1 \over \mu}\psi(0,\cdot) + ({1 \over \mu} - {g \over \lambda})(\varphi +
\varphi^t)(0)\bigr)\cr &\quad \quad + \bigl( {1 \over \lambda}\psi(\cdot,0) - {f
\over\mu}(\varphi +
\varphi^t)(0)
\bigr) \otimes E_+\cr
& \quad \quad +  {f \over \mu}\sum_{i \neq +} [E_+, E_{-i}] \otimes E_i  
  + {g \over \lambda} \sum_{i \neq +} E_{i}
\otimes [E_+,E_{-i}].\cr
}\eqno(9.14)
$$ 
This yields $g = xf$ and the result is that
$$
 r = \varphi + \sum E_{-i} \otimes E_i + {x \over 1-x} \, C,\quad x = \lambda/\mu,
\eqno(9.15)
$$ 
which agrees with the simplest r-matrix  in [BD], but in the notation of [J].  

\ve

\no {\bf 10. The Standard Classical r-matrix for Twisted 
Loop Algebras.}

The construction of a twisted affine Kac-Moody algebra [K] involves two simple Lie algebras,
\g ~and a subalgebra \g $_0$, such that \g ~admits a diagram automorphism
of order $k = 2$ or $3$  to which is associated   a Lie algebra automorphism $\mu$ that
centralizes  \g $_0$. The eigenvalues of $\mu$ are of the form $\omega^j,\, j =
0,1,\ldots $ , and \g $ = \sum_{j=0}^{k-1} $\g$_{j}$, where 
\g$_{j}$ is the
sum of the eigenspaces  with eigenvalues $ \omega^{j\,{\rm mod}\,k}$.
 The restriction of the adjoint action of
\g to \g$_0$ acts irreducibly on each \g $_j$.

Now let $\{H_a, e_{\pm \alpha}\} \alpha = 1,\ldots n$ be a Chevalley basis for 
\g$_0$, and let $E_+$ be a highest weight vector (for the action of \g$_0$) in
\g$_1$. Then $\{e_\alpha\},E_-$ generate \g, and
$$ [e_\alpha,E_+] = 0 = [e_{-\alpha},E_-].\eqno(10.1)
$$
 The twisted loop algebra~ $ \hat {\hbox{\g}} = \Crm [\lambda, {1 \over \lambda}] \,\otimes\,$ \g
 ~  is generated by
$\{e_{\pm \alpha}\}, \alpha = 0,\ldots ,n$, with
$$ 
e_0 = \lambda E_-,\quad e_{-0} = {1 \over \lambda}E_+.
\eqno(10.2)$$
This algebra is of the type ${\cal A}_{cl}'$, so our standard R-matrix applies. We define $r$
in terms of the expansion of $R$ in powers of $\hbar$ and work out the implications of
the relations (9.7).

Let $\{E_i\}$ be a Weyl basis for \g$_0$ and normalize so that the Casimir element
for that algebra is
$$
C_0 = \varphi + \varphi^t + \sum E_{-i} \otimes E_i + \sum E_i \otimes E_{-i}.\eqno(10.3)
$$
Then a special solution of (9.7) with $\gamma \neq 0$ is given by the first two terms in
(9.13) and the general solution is
$$
r = \varphi + \sum E_{-i} \otimes E_i + \sum_0^{k-1} f_j C_j,
$$
where $C_j$ is the projection of the Casimir element $C$ of \g ~on \g$_j$,
on the first factor. Now (9.7), with $\gamma = 0$:
$$
\eqalign{
0 &= [ 1 \otimes e_{-0} + e_{-0} \otimes 1, \, \sum E_{-i} \otimes E_i + 
\sum f_j C_j] + \sum e_{-\alpha} \otimes [e_\alpha, e_{-0}]\cr
&= { 1 \over \mu}\sum [E_+,E_{-i}] \otimes E_i + \sum f_j \bigl(
{1\over \lambda} [1 \otimes E_+, C_j] + {1 \over \mu}[E_+ \otimes 1, C_j]\bigr)\cr
&\hskip1in+ {1 \over \mu} E_+ \otimes (\varphi + \varphi^t)(0)\cr
&= {1 \over \mu}[E_+ \otimes 1,C_o] + \sum \bigl({f_j \over \lambda} [ 1\otimes E_+, C_j] +
{f_{j-1}\over \mu} [ E_+ \otimes 1, C_{j-1}].
\cr}\eqno(10.4)
$$
This vanishes iff
$$
\eqalign{
& f_1 = x(f_0 +1),\,\, f_0 = x f_1,\,\,\, k = 2,\cr
& f_1 = x(f_0 +1),\,\, f_2 = x f_1, \,\, f_0 = x f_2,\,\,\, k = 3,\cr
}
$$
That is,
$$
f_j = {x^j \over 1 - x^k } \, C_j - \delta_j^0 \, C_0.
$$
Finally, the unique solution is 
 $$
r = \varphi + \sum E_{-i} \otimes E_i - C_0 + {1 \over 1-x^k}\sum_0^{k-1} x^j C_j,\eqno(10.5)
$$
again in agreement with [BD], in the notation of [J].
\vskip.50cm
\no {\bf Remark.} Choose a basis of weight vectors in \g$_1$, then 
$$
C_1 = E_- \otimes E_+ + E_+ \otimes E_- + \ldots ,
$$
with unit coefficients for the contributions with highest weight. This follows from the
normalization in (10.3) and fact that $1 \otimes E_+ +  E_+ \otimes  1 $ commutes
with $C = \sum C_j$.

\vskip1.5cm

\no {\bf 11. Including the central extension.}

\underbar {The untwisted case.} The extension is recovered by omitting the replacement of
$\hat\varphi$ by $\varphi$ in (9.3). We can still represent the r-matrix as a power series in $x
= \lambda/\mu$, but it is no longer true, as it was in the case of the loop group, that
$[e_0,e_{-0}] = [E_-,E_+]$. Instead,
$$
[e_0,e_{-0}] = (\hat \varphi + \hat \varphi^t)(0) = [E_-,E_+] + c.\eqno(11.1)
$$
More generally, for polynomials $f,g \in \Cit[\lambda,{1 \over \lambda}]$, and $x,y \in
{\cal A}'_{cl}$,
$$
[fx,gx] = fg[x,y] + c\,<x,y>\, {\rm Res} (f'g),\eqno(11.2)
$$
where the form $<,>$ is the invariant form on ${\cal A}'_{cl}$ normalized as follows:
If the Casimir element is $C^{ij} x_i \otimes x_j$, then $<x_i,x_j> = (C^{-1})_{ij}$;
Res$(f)$ is the constant term in $\lambda f$.
\vskip.50cm

\no {\bf Remark.} This normalization implies that
$$
[fC_{12},gC_{23}] = fg [C_{12},C_{23}] + c_2 C_{13} {\rm Res} (f'g).\eqno(11.3)
$$ 

\vskip.50cm

This change leaves (9.12) and (9.13) unaffected, while (9.14) becomes 
$$
\eqalign{
0 &= E_+ \otimes \bigl({1 \over \mu}\psi(0,\cdot) + {1 \over \mu} (\hat \varphi +
\hat \varphi^t)(0) + [e_0,g(x)E_-]\bigr)\cr 
&\quad \quad + \bigl( {1 \over\lambda}\psi(\cdot,0) - {f
\over\mu}(\varphi +
\varphi^t)(0)
\bigr) \otimes E_+\cr
& \quad \quad +  {f \over \mu}\sum_{i \neq +} [E_+, E_{-i}] \otimes E_i  
  + {g \over \lambda} \sum_{i \neq +} E_{i}
\otimes [E_+,E_{-i}].\cr
} 
$$ 
The modification in the second term ($\varphi$ replaced by $\hat \varphi$) is exactly
compensated by a new contribution from the linear $\lambda$-term in $g$. (There 
is no linear $\mu$-term in $f$.) The conclusion is that the new r-matrix is
$$
 \hat r = \hat \varphi + \sum E_{-i} \otimes E_i + {x \over 1-x} \, C. \eqno(11.4)
$$ 
\no \underbar {The twisted case.} It is easy to verify, with the help of the remark at the
end of Section 3, that the restitution $\varphi \rightarrow \hat \varphi$ can be made without
affecting the cancellations; so the result is that
 $$
\hat r = \hat \varphi + \sum E_{-i} \otimes E_i - C_0 + {1 \over 1-x^k}\sum x^j C_j.\eqno(11.5)
$$

It is amusing to verify directly that the classical Yang Baxter relation for $r$,
$$
\YB(r) := [r_{12}, r_{13} + r_{23}] + [r_{13}, r_{23}] = 0,
$$
implies the same relation for $\hat r$: The inclusion of the extra term in $\hat \varphi$ means
that
$$
\YB(\hat r) = \YB(r) + [r_{13},(c \otimes d)_{23}].\eqno(11.6)
$$
The evaluation of $\YB(r)$ now has to take into account the new term (involving $c$) in 
Eq.(11.2). Actually, only $[r_{12},r_{23}]$ is affected, and with the aid of Eq.(11.3) one
finds that the new contribution is
$$
\YB(r) = c_2 \lambda {d \over d\lambda} r_{13},
$$ 
which exactly cancels the other term. In the twisted case one must use the following
generalization of Eq.(11.3):
$$
[f{C_j}_{12},g{C_{j'}}_{23}] = fg [{C_j}_{12},{C_{j'}}_{23}] + \delta_j^{j'}c_2
\,{C_j}_{13} {\rm Res} (f'g).\eqno(11.7)
$$

\ve

 \font\head=cmss8 scaled\magstep1
{\ce{\head DEFORMATIONS}}
\vskip.5cm

\line {\bf 12. First Order Deformations. \hfil}
 
Quantum groups can be understood as deformations of the Hopf structure associated with
Lie algebras or Kac-Moody algebras. The point of view that emphasizes the direct  connection 
between quantum groups and Lie groups, as well as the deep roots of quantum groups in 
deformation theory and in the theory of $*$-products, has been shown to lead to profound 
insight into their general structure [BFGP][BP][EK]. Here we use deformation
theory with a  different purpose. The initial structure is the bialgebra
associated with a standard R-matrix, with a fixed set of parameters. The
deformed structure is a bialgebra equipped with an R-matrix that is non-standard
and that does not commute with the Cartan sub-algebra. We emphasize that the
context is more general than quantized Kac-Moody algebras.

This work was initiated with the aim of calculating the universal
R-matrices associated with simple Lie algebras, as deformations of the
standard universal R-matrix.  We shall establish a direct correspondence
between the classical r-matrices of Belavin and Drinfeld on the one hand,
and the deformations of the standard, universal R-matrix for simple
quantum groups on the other. In preparation for this we have explored the
meaning of the Yang-Baxter relation in a much more general context, and
we shall endeavor to maintain this generality in our approach to
deformations.  But, as for the types of deformations,   we shall   limit
our study in a way that seems natural in the context of quantum groups.

A deformation of the standard R-matrix is a formal series
$$ 
R_\epsilon=R+\epsilon R_1+\epsilon^2 R_2 + \ldots~. \eqno(12.1)
$$
\no Here $R$ is a standard R-matrix on ${\cal A}' = {\cal A}/I$ with any choice of 
parameters and the ideal $I$ determined by them. The coefficients  $t_{(\alpha)}^{(\alpha')}$ of $R$
are determined by the Yang-Baxter relation, and we attempt to find
$R_1,R_2,\ldots$ so that $R_\epsilon$ will satisfy the same relation to
each order in
$\epsilon$.  To make this program precise, we must specify the nature of
the leading term; the remainder should then be more or less unique.
\ve
Recall that $R$ ``commutes with Cartan."  An element $Q\in
{\cal{A}}\otimes {\cal{A}}$ is said to have weight $w$ if
$$ [H_a\otimes 1+1\otimes H_a,~Q\,] = w_a Q~, \quad w_a\in \Crm,~a\in M~.
\eqno(12.2)
$$
\no Thus $R$ has weight zero.  The image of $Q$ by the projection
${\cal{A}}\otimes {\cal{A}}\rightarrow {\cal{A}}^\prime\otimes
{\cal{A}}^\prime$ has the same weight.  We shall suppose that $R_1$ is
homogeneous (has weight), but this restriction is inessential and will be relaxed later.

Recall further that $R$ is driven by the linear term; by virtue of the
Yang-Baxter relation, $R$ is completely determined by the term
$e_{-\alpha}\otimes e_\alpha$.  It is natural to study deformations that
are driven by a similar term, with fixed, non-zero weight:
$$
 R_1=S(e_{\pm\sigma} \otimes e_{\pm\rho}) + \ldots~, \eqno(12.3)
$$
\no with $\sigma,\rho$ fixed and the factor $S$ is in ${\cal{A}}^0$. (The
unwritten terms are of higher order, in a sense that we shall make precise in a moment.) Such deformations may be called  
``non-singular", to contrast them to singular deformations for which the term of order $\epsilon$ is either
absent or else of a form that sets it appart from the driving term in the undeformed R-matrix. 
We do not claim that that this exhausts the possibilities.
In fact, we know of a ``singular" deformation that is driven by an $R_1$ of higher order in the generators [FG2].
It is highly special and occurs only when some of the parameters are roots of unity. But we 
believe that the deformations studied here have the best chance of possessing a 
cohomological interpretation.

We shall now make precise the concept of ``higher order". 

\vskip.5cm

\no {\bf Proposition 12.}  The algebra ${\cal{A}}^\prime={\cal{A}}/I$ is
Z\hskip-1.5mm Z-graded, with grade $e_{\pm\alpha}=\pm 1$, grade
$H_a=0$.  An alternative grading is obtained by reversing the sign.
\vskip.5cm

{\bf Proof.} This is a consequence of the fact that the generators of $I$
are homogeneous; ${\cal A}'$ inherits the grading of $\cal A$. 

\vskip.5cm

The standard R-matrix is a formal series $\sum_k \psi^-_k\otimes
\psi^+_k,~\psi^\pm_k\in {\cal{A}}^\prime$.  We use the grading of
Proposition 12 in the second space, the alternative grading in the first
space; then grade $\psi^\pm_k=k$ and $R$ is a formal sum of  terms with
grade $(k,k),~k=0,1,2,\ldots~.$  This grading is an extension of that
used previously, made necessary by the appearance of $e_\sigma$ in the
first space and $e_{-\rho}$ in the second.

With the inclusion of (12.3)   the grades  descend to
$(-1,-1)$.  Finally, the unwritten terms in (12.3) is a series by
ascending grades.  The fact that the grades are bounded below is
fundamental.  We claim that $R_\epsilon$, a formal series in
$\epsilon$, each term a formal series in ascending grades, if it
satisfies the Yang-Baxter relation, is completely determined by the
choice of the two generators $e_{\pm\sigma}$ and $e_{\pm\rho}$ in (12.3).

We shall see that the standard R-matrix on ${\cal A}'$, with the
parameters of 
${\cal{A}}'$ in general position, is rigid with respect to deformations
of the type (12.3).  We begin our investigation by establishing some
conditions on the parameters that are necessary for the existence of a
deformation.  We shall study each of the
four possibilities envisaged by (12.3) separately.  We organize the
contributions to
$$ Y\hskip-1.0mm B_\epsilon~:=~ R_{\epsilon 12}R_{\epsilon 13}R_{\epsilon
23}- R_{\epsilon 23}R_{\epsilon 13}R_{\epsilon 12}
$$
\no in the same way as the contributions to $Y\hskip-1.0mm B$.  A term
$\psi_1\otimes \psi_2\otimes\psi_3$ is said to have grade $(\ell,n)$ if
$\psi_3$ has grade $n$ and $\psi_1$ has alternative grade
$\ell$. We limit ourselves  to terms linear in $\epsilon$ and end this section by
disposing of three of the four possibilities in (12.3).
\vskip.5cm

 \no \underbar  {Deformations of Types} $e_{-\sigma}\otimes e_\rho$,
$e_\sigma \otimes e_\rho$ \underbar {and} $e_{-\sigma} \otimes e_{-\rho}$.  
~Suppose first that the driving term in $R_1$ is
$$ S(e_{-\sigma}\otimes e_\rho)~, \quad S\in {\cal{A}}^0\otimes
{\cal{A}}^0~.
$$
\no We examine the contributions to $  Y\hskip-1.0mm B_\epsilon$ of order
$\epsilon$.

The lowest grades are (1,0) and (0,1), with contributions
$$
\eqalign{&(Se_{-\sigma}\otimes e_\rho)_{12}~
R^0_{13}R^0_{23}-R^0_{23}R^0_{13}(Se_{-\sigma}\otimes e_\rho)_{12}~, \cr
&R^0_{12}R^0_{13}(Se_{-\sigma}\otimes e_\rho)_{23}- (Se_{-\sigma}\otimes
e_\rho)_{23}R^0_{13}R^0_{12}~, \cr}
$$
\no respectively.  These vanish if and only if
$$
e^{\varphi(\sigma,\cdot)-\varphi(\rho,\cdot)}=1=e^{\varphi(\cdot,\sigma)-\varphi(\cdot,\rho)}~.
\eqno(12.4)
$$
\no In grade (1,1) we encounter additional restrictions, 
$$ 
e^{\varphi(\rho,\cdot)+\varphi(\cdot,\sigma)}=1~. \eqno(12.5)
$$
 Conditions (12.4)-(12.5) are necessary.   It follows   that
$$
e^{\varphi(\rho,\alpha)+\varphi(\alpha,\rho)}=1=e^{\varphi(\sigma,\alpha)+\varphi(\alpha,\sigma)}~,
\quad \alpha\in N~.  
$$
\no These are conditions that are familiar from our investigation of
constants, see Eq.(3.8).  The relations that are thus implied are
$$ e_\rho e_\alpha-e^{\varphi(\rho,\alpha)}e_\alpha e_\rho = 0 = e_\sigma
e_\alpha -e^{\varphi(\sigma,\alpha)}e_\alpha e_\sigma~,
\,\, \forall\alpha \in N.
$$
  This constitutes a high degree of commutativity in ${\cal{A}}^\prime$
and takes us far away from our main interest in simple quantum groups. 
We therefore end our investigation of the type $e_{-\sigma}\otimes
e_\rho$ at this point.
 
  Similar results are obtained for
deformations of type
$e_{-\sigma}\otimes e_{-\rho}$ and $e_\sigma \otimes e\rho$. 

\vskip1.5cm

\line {{\bf 13. First Order Deformations of Type $e_\sigma\otimes e_{-\rho}$ 
and the Classical Limit.} \hfil}
\s

We come to the last case envisaged in Section 12.  Eq. (12.3), when the
driving term in $R_1$ has the form
$$ S(e_\sigma\otimes e_{-\rho})~, \quad S\in {\cal{A}}^0\otimes
{\cal{A}}^0~. \eqno(13.1)
$$
\no This term has grade (-1,-1); it is the only term in $R_1$ with this
grade, the lowest.  The factor $S$, and all other terms in
$R_1$, are completely determined by the Yang-Baxter relation
$  Y\hskip-1.0mm B_\epsilon=0$ to first order in $\epsilon$.  Besides
(13.1) there is in $R_1$ one other term with only two roots, of the form
$$ S^\prime(e_{-\rho}\otimes e_\sigma)~, S^\prime\in {\cal{A}}^0\otimes
{\cal{A}}^0~;
\eqno(13.2)
$$
\no it has grade (1,1).
\vskip.5cm

\ve

\no {\bf Theorem 13.1.}  Let $R$ be the standard R-matrix described in
Theorem 7.  Suppose that $R+\epsilon R_1$ is a first order deformation,
satisfying the Yang-Baxter relation to first order in
$\epsilon$.  Suppose also that the term of lowest grade in $R_1$ has
the form (13.1); then  the parameters satisfy
$$ 
e^{\varphi(\cdot,\rho)+\varphi(\sigma,\cdot)}=1~. \eqno(13.3)
$$ 
Conversely, when the parameters are in general position on this surface, then there exists a unique first 
order deformation such that the term of lowest grade has the form (13.1), namely
$$ 
R_1=  R(Ke_\sigma\otimes Ke_{-\rho}) - (Ke_{-\rho}\otimes Ke_\sigma)R
~, \eqno(13.4)
$$
\no with $K~:=~e^{\varphi(\cdot,\rho)}$.
\b

\no {\bf Proof.}  An easy calculation in the lowest grades shows that
(13.3) is necessary and that $S=K\otimes K$, up to a  numerical factor
that we fix once and for all.

Let $R_1^i,~i=1,2,$ be the two summands in (13.4).  The term of order
$\epsilon$ in $ Y\hskip-1.0mm B_\epsilon$ is the sum of the following six
quantities:
$$
\eqalign{A^i_{12} &= (R^i_1)_{12}R_{13}R_{23}- R_{23}R_{13}(R^i_1)_{12}~,
\cr A^i_{13} &= R_{12}(R^i_1)_{13}R_{23}-R_{23}(R^i_1)_{13}R_{23}~, \cr
A^i_{23} &= R_{12}R_{13}(R^i_1)_{23}-(R^i)_{23}R_{13}R_{12}~, \quad
i=1,2~.
\cr}
\eqno(13.5)
$$
\vskip.5cm
\no {\bf Step 1.}  We begin with the term that contains the lowest grade,
(1,1):
$$
\eqalign{&A^1_{13}=R^i[-\alpha]R^j[-\beta]Ke_\sigma\otimes R_i[\alpha]
R^k[-\gamma]\otimes R_j[\beta]Ke_{-\rho}R_k[\gamma]-\ldots~,
\cr &[-\alpha]\otimes [\alpha]~:=~e_{-\alpha_1}\ldots e_{-\alpha_\ell}
\,t^{(\alpha^\prime)}_{(\alpha)}\,e_{\alpha^\prime_1}\ldots
e_{\alpha^\prime_\ell}~. \cr}
\eqno(13.6)
$$
\no A sum over indices and numbers of indices ($\ell~\alpha$'s, 
$m~\beta$'s and $n~\gamma$'s) is understood, and $-\ldots$ stands for the
reflected term.  Using the fact that $R$ satisfies $  Y\hskip-1.0mm B=0$
we can convert (13.6) to
$$ A^1_{13}=R^i[-\alpha]R^j[-\beta]Ke_\sigma\otimes R_i[\alpha]
R^kK[-\gamma]\otimes R_j[\beta]KR_k\bigl[e_{-\rho},[\gamma]\bigr]~ +
\ldots ,
\eqno(13.7)
$$
\no where $+\ldots$ stands for a similar expression that contains a
factor $\bigl[e_\sigma,[-\alpha]\bigr]$ in the first space.  We have used
(13.3) and continue to use this relation without comment.
\vskip.5cm

\no {\bf Step 2.}  Evaluate the commutators in (13.7) using (5.1) and
(5.7).  The result
$$
\eqalign{A^1_{13}&=R^i[-\alpha]R^j[-\beta]Ke_\sigma\otimes
R_i[\alpha]R^kK[-\gamma]e_{-\rho}\otimes R_j[\beta]KR_k[\gamma]K^{-1}  +
\ldots \cr} \eqno(13.8)
$$
\no is a sum of four similar expressions.  Note that the evaluation of
the commutators involves a shift in the summation indices
$\ell,m,n$. The generators $e_\sigma, e_{-\rho}$, in spaces 1,3 in (13.7), are now in spaces 1 and 2,
and the lowest grades in (13.8) are (-1,0) and (0,-1).
\vskip.5cm

\no {\bf Step 3.}  Now write down the full expression for
$A^1_{12}+A^1_{23}$; it also contains four similar terms.  Two of them
cancel two of the terms in (13.8), by virtue of the relation  $Y\hskip-1.0mm B=0$.
\vskip.5cm

\no {\bf Step 4.}  Combine the remaining two terms from (13.8) with  the
remaining two terms from $A^1_{12}+A^1_{23}$ and verify that
$$
\eqalign{&A^1_{12}+A^1_{13}+A^1_{23}  \cr &=
R^i[-\alpha]KR^j\bigl[e_\sigma,[-\beta]\bigr]\otimes R_i[\alpha]
KR^ke_{-\rho}[-\gamma]\otimes
R_jK[\beta]R_ke^{\varphi(\rho,\cdot)}[\gamma]  + \ldots \, ,\cr}
\eqno(13.9)
$$
\no where $+\ldots$ stands for a term that contains a factor
$\bigl[[\beta],e_{-\rho}\bigr]$ in the third space.
\vskip.5cm

\no {\bf Step 5.}  Evaluate the commutators (second shift of summation
indices)
$$
\eqalign{&= R^i[-\alpha]R^jK[-\beta]K^{-1}\otimes
R^i[\alpha]KR^ke_{-\rho}[-\gamma]\otimes R_jK[\beta]e_\sigma
R_ke^{\varphi(\rho,\cdot)}[\gamma] + \ldots \cr &=:~X_1+X_2+Y_1+Y_2~.
\cr}  \eqno(13.10)
$$ The lowest grades are now (1,0) and (0,1). The generator $e_\sigma$ has completed its journey
towards the east and is found in the third space; the generator $e_{-\rho}$, travelling westward, 
is in space two.
\vskip.5cm

\no {\bf Step 6.}  Two of the four terms in (13.10) are:
$$
\eqalign{ X_1 &= R_{12}R_{13}\bigl\{(Ke_{-\rho}
\otimes Ke_\sigma)R\bigr\}_{23}~, \cr Y_2 &=
-R_{23}R_{13}\bigl\{(Ke_{-\rho}\otimes Ke_\sigma)R\bigr\}_{12}~. \cr} 
\eqno(13.11)
$$
\no Now add $A^2_{12}+A^2_{23}$ to (13.10) to get
$$
\eqalign{A^1_{12}&+A^1_{13}+A^1_{23}+A^2_{12}+A^2_{23} \cr &= \tilde
X_1+X_2+Y_1+\tilde Y_2~, \cr} \eqno(13.12)
$$
\no where $\tilde X_1$ and $\tilde Y_2$ are obtained from $X_1$ and
$Y_2$ by adding $A^2_{23}$ and $A^2_{12}$.
\ve

\no {\bf Step 7.}  Use the relation $Y\hskip-1.0mm B=0$ to modify the
expressions for $\tilde X_1$ and $\tilde Y_2$; then notice that the four
terms in (13.12) can be combined to two,
$$ = R^iK[-\alpha]R^j[-\beta] \otimes
R_iK[e_{-\rho},[\alpha]]R^k[-\gamma] \otimes Ke_\sigma
R_j[\beta]R_k[\gamma] + \ldots,
\eqno(13.13)
$$ where the other term has a factor $[[-\gamma],e_\sigma]$ in the second
space.
\vskip.5cm

{\bf Step 8.} Evaluate the commutators (third shift of summation indices),
$$ = R^iK[-\alpha]e_{-\rho}R^j[-\beta] \otimes
R_jK[\alpha]K^{-1}R^k[-\gamma] \otimes K e_\sigma R_j[\beta]R_k[\gamma] +
\ldots.
\eqno(13.14)
$$ This expression has four terms; the lowest grade is (1.1). The generators $e_\sigma,e_{-\rho}$ 
have reached their final destination, $e_\sigma$ is in space three and $e_{-\rho}$ is in space one.
What remains can   be compared with the last of the six contributions to $Y\hskip-.8mm B_\epsilon$,
namely $A^2_{13}$.
\vskip.5cm

{\bf Steps 9, 10.} Two of the four terms in (13.14) cancel each other because
$Y\hskip-1.0mm B = 0$ and
 the remaining two terms add up to $-A_{13}^2$.
\vskip.5cm

This completes the verification of the claim that (13.4) defines a first
order deformation of $R$. To complete the proof of Theorem 13.1 we must
show that this expression (13.4) is unique. This was done by complete
mathematical induction. We omit the details but point out that the key to
the induction processs is visible in steps 2,5 and 8, where the summation
indices are shifted. Theorem 13.1 is proved.
\vskip.5cm

Let ${\cal P}$ be the collection of pairs $(\sigma,\rho) \in N \otimes N$
such that (13.3) holds; each distinct pair defines a first order
deformation $R + \epsilon R_1^{\sigma,\rho}$ of $R$. Because these
deformations are only first order they generate a linear space 
$$ R_1 = \sum_{\sigma,\rho \in {\cal P}} C_{\sigma,\rho}
R_1^{\sigma,\rho},
\eqno(13.15)
$$  with coefficients in \Crm. The dimension of this space of first
order deformations is zero for parameters in general position. It remains
zero, generically, when the parameters are such that the ideal $I$
generated by the constants is non-zero and
$R$ is  defined on ${\cal A}/I$. The exceptional points in the space of parameters,
at which there are pairs $(\sigma,\rho)$ satisfying (13.3), are bifurcation points in 
the space of generalized quantum groups.

To any first order deformation of $R$, there corresponds a first order
deformation of $r$,\break
$$R_{\epsilon} = 1 + \hbar r_\epsilon + o(\hbar^2), \quad 
r_\epsilon=r+\epsilon r_1+o(\epsilon^2)~. \eqno(13.16)
$$
\no Eqs. (13.4) and (13.15) give us
$$ r_1=\sum_{\sigma,\rho\in {\cal{P}}} C_{\sigma,\rho} (e_\sigma \wedge
e_{-\rho})~, \eqno(13.17)
$$
\no where ${\cal{P}}$ is the set of pairs with the property
$$
\varphi(\rho,\cdot)+\varphi(\cdot,\sigma)~:=~0~. \eqno(13.18)
$$

The original work of Belavin and Drinfeld culminates in a list of
constant r-matrices that is complete up to equivalence.  Their results
have recently been re-derived in terms of deformation theory and the
associated cohomology.
\vskip.5cm

\no {\bf Proposition 13.} [F] Let $r$ be the standard r-matrix (8.14)
for a simple Lie algebra ${\cal L}$.  The space of essential, first
order deformations of $r$, satisfying (8.15) and (8.16), is
$$ H^2({\cal L}^*,\Crm)=\bigl\{r_1=\sum_{\sigma,\rho\in {\cal{P}}}
C_{\sigma,\rho}e_\sigma \wedge e_{-\rho}+ \sum \tilde C^{ab} H_a\otimes
H_b\bigr\}~. \eqno(13.19)
$$
\no The exact deformations are of finite order and coincide with the
r-matrices of [BD].

\vskip.5cm The second, Cartan term is not ``essential" in the present
context; it represents the freedom to vary the parameters.  We conclude
that
\vskip.5cm

\no {\bf Theorem 13.2.}  The first order deformations of the standard
R-matrix described in Theorem 13.1, upon specialization to a simple
quantum group, are in one-to-one correspondence, via (13.16), with the
  first order essential deformations of the associated standard r-matrix,
modulo variations of the parameters.
\vskip.5cm
One concludes that the class of deformations investigated in Section 13 is wide enough to 
encompass the quantization of all simple Lie bialgebras. We shall see that the affine 
Kac-Moody algebras are provided for also.

\vskip1.5cm

\line {\bf 14. Hopf Structure. \hfil}
\s

It is of some interest to verify that the standard R-matrix, satisfying
the Yang-Baxter relation, actually intertwines the coproduct of a Hopf
algebra with its opposite.
\vskip.5cm

\no {\bf Proposition 14.1.}  (a) There exists a unique homomorphism
$\Delta:~{\cal{A}}\to {\cal{A}}\otimes {\cal{A}}$, such that
$$
\eqalign{\Delta(H_a) &= H_a\otimes 1+1\otimes H_a~, \quad a\in M~, \cr
\Delta(e_\alpha) &= 1\otimes e_\alpha + e_\alpha \otimes 
e^{\varphi(\alpha,\cdot)}~, \cr
\Delta(e_{-\alpha}) &= e^{-\varphi(\cdot,\alpha)} \otimes
e_{-\alpha}+e_{-\alpha}\otimes 1~, \quad \alpha \in N~. \cr}
\eqno(14.1)
$$
\no (b) If $I\subset {\cal{A}}$ is the ideal generated by the constants
in ${\cal{A}}^+$ and ${\cal{A}}^-$, and ${\cal{A}}^\prime= {\cal{A}}/I$,
then $\Delta$ induces a unique homomorphism
${\cal{A}}^\prime \to {\cal{A}}^\prime\otimes {\cal{A}}^\prime$ that will
also be denoted $\Delta$, so that (14.1) holds with $H_a$ and
$e_{\pm\alpha}$ being interpreted as generators of ${\cal{A}}/I$.
\s
\no (c) Let $\Delta^\prime$ be the opposite coproduct on ${\cal{A}}/I$,
and $R$ the standard R-matrix on ${\cal{A}}/I$ (satisfying Yang-Baxter),
then $\Delta R=R\Delta^\prime$.
\s
\no (d) The algebra ${\cal{A}}$ becomes a Hopf algebra when endowed with
the counit ${\cal E}$ and the antipode $S$.  The former is the unique
homomorphism ${\cal{A}}\to \Crm$ such that
$$
\eqalign{& {\cal E}(a)=1~, \quad {\cal E} (H_a)=0~, \quad a\in M~, \cr &
{\cal E}(e_{\pm\alpha})=0~, \quad \alpha\in N~. \cr}
\eqno(14.2)
$$
\no The antipode is the unique anti-automorphism $S:~{\cal{A}}\to 
{\cal{A}}$ such that
$$
\eqalign{& S(1)=1~, \quad S(H_a)=-H_a~, \quad a\in M~, \cr &
S(e_\alpha)=-e_\alpha e^{-\varphi(\alpha,\cdot)}~, \quad
S(e_{-\alpha})=-e^{\varphi(\cdot,\alpha)} e_{-\alpha}~, \quad
\alpha\in N~. \cr} \eqno(14.3)
$$
\no (e) The counit ${\cal E}$ and the antipode $S$ of ${\cal{A}}$ induce
analogous structures on ${\cal{A}}^\prime={\cal{A}}/I$ such that (14.3)
holds on ${\cal{A}}^\prime$.
\vskip.5cm

\no {\bf Proof. }
 (a) The verification amounts to checking that $\Delta({\cal{A}})$ has
the relations of ${\cal{A}}$, in particular,
$$
\bigl[\Delta(e_\alpha),\Delta(e_{-\beta})\bigr] =
\delta^\beta_\alpha~\Delta\bigl([e_\alpha,e_{-\beta}]\bigr)~.
$$ (b) The ideal $I$ is generated by elements $x\in {\cal{A}}^+$ and
$y\in {\cal{A}}^-$ such that 
$[e_{-\alpha},x]=0=[e_\alpha,y],~\alpha\in N$.  Since 
$\Delta:~{\cal{A}}\to {\cal{A}}\otimes {\cal{A}}$ is a homomorphism,
$\Delta$ induces a homomorphism ${\cal{A}}/I\to ({\cal{A}}\otimes
{\cal{A}})/\Delta(I)$.  We must show that 
$\Delta(I)\subset I\otimes {\cal{A}} + {\cal{A}}\otimes I$.  Since $I$ is
generated by elementary constants, it is enough to show that, for an
elementary constant $C$, $\Delta (C) \subset I\otimes {\cal{A}} +
{\cal{A}}\otimes I$.  Let $C\in {\cal{A}}^+$ be an elementary constant;
then $[e_{-\alpha},C]=0$ and thus $[\Delta(e_{-\alpha}),\Delta (C)]=0,~
\alpha\in N$.  If $C$ is of order $n$ in the generators,  (14.1)
shows that
$$
\Delta C=1\otimes C+ P^1\otimes P_{n-1}+P^2\otimes P_{n-2}+\ldots +
C\otimes P_0~,
$$
\no where $P^n$ and $P_k$ are homogeneous of order $k$ in the
$e_\alpha$'s.  Because $C$ is an elementary constant---Definition 4.1.---
there is no constant among the $P^k,P_k,~n=1,\ldots,n-1$; then
$[\Delta (e_{-\alpha}),\Delta (C)]=0$ implies that $\Delta C= 1\otimes
C+C\otimes P_0$ which indeed belongs to $ I\otimes {\cal{A}} +
{\cal{A}}\otimes I$; consequently $\Delta$ provides a map ${\cal{A}}/I\to
{\cal{A}}/I\otimes {\cal{A}}/I$.

(c) We use the abbreviation~--~ compare (13.6), Definition 2.2 and Eq.
(2.9)~-
$$ R=t_{(\alpha)}^{(\alpha^\prime)} R^i[e_{-\alpha}]\otimes 
R_i[e_{\alpha^\prime}]~,
$$
$$
\eqalign{\Delta (e_\beta)R-R\Delta^\prime (e_\beta) &=
t^{(\alpha^\prime)}_{(\alpha)}
\bigl(R^i[e_{-\alpha}]\otimes e_\beta R_i[e_{\alpha^\prime}] \cr
&+e_\beta R^i[e_{-\alpha}]\otimes e^{\varphi(\beta,\cdot)}
R_i[e_{\alpha^\prime}]-R^i[e_{-\alpha}]e_\beta\otimes
R_i[e_{\alpha^\prime}] \cr
&-R^i[e_{-\alpha}]e^{\varphi(\beta,\cdot)}\otimes R_i[e_{\alpha^\prime}]
e_\beta\bigr)~. \cr}
$$
\no Terms 2 and 3 combine to $R^i[e_\beta,t^{(\alpha^\prime)}]\otimes
R_i[e_{\alpha^\prime}]$, and the recursion relations (7.2) implies that
the sum of all four terms equals zero. Actually this recursion relation, when summed over $n$,
is nothing more than the statement $\Delta (e_\beta)R-R\Delta^\prime (e_\beta) = 0$. It should 
 be pointed out that the co-product was not known {\it a priori}; the Yang-Baxter relation 
gave us the recursion relation and this amounts to a   determination of the co-product.

(d) The existence and   uniqueness of the homomorphism ${\cal E}$ and
the anti-homomor-phism $S$ are obvious.  We have to show that
${\cal E}$ satisfies the axioms
$$ ({\cal E}\times id)\Delta=id=(id\times{\cal E})\Delta~,
$$
\no which is straightforward, and that
$$ m(id\times S)\Delta=\epsilon=m(S\times id)~.
$$
\no Here $m$ indicates multiplication, ${\cal{A}}\otimes {\cal{A}}\to
{\cal{A}}$.  For example,
$$
 m(id\times S)\Delta (e_\alpha)=S(e_\alpha)+e_\alpha
e^{-\varphi(\alpha,\cdot)}=0~.
$$ (e) Obvious, since ${\cal E}(I)=0$ and $S(I)=I$ by Proposition 7.2.
Proposition 14.1 is proved.
\ve
We turn to the case of the deformed R-matrix of Section 13, all
statements should be understood to hold to first order in the deformation
parameter $\epsilon$.  The maps $\Delta,{\cal E}$ and $S$ are    as
before and the deformed maps are
$
\Delta_\epsilon=\Delta + \epsilon\Delta_1~, \quad {\cal E}_\epsilon={\cal
E}+{\epsilon}{\cal E}_1~, \quad S_\epsilon=S+\epsilon S_1~.
$
\vskip.5cm
\no {\bf Proposition 14.2.}  (a) There is a unique homomorphism
$\Delta_\epsilon:~{\cal{A}}\to {\cal{A}}\otimes {\cal{A}}$ such that
$$
\Delta_1(x)=[\Delta(x),~Ke_{-\rho}\otimes Ke_\sigma]~, \quad x\in
{\cal{A}}~.\eqno(14.4)
$$
\no (b) The projection of $\Delta_\epsilon$ to ${\cal{A}}^\prime\to
{\cal{A}}^\prime \otimes {\cal{A}}^\prime$ is well defined. (c) Let
$\Delta^\prime_\epsilon$ be the opposite coproduct on
${\cal{A}}^\prime={\cal{A}}/I$, and $R_\epsilon = R+\epsilon R_1$ the
R-matrix of Theorem 13.1, then $\Delta_\epsilon R_\epsilon= R_\epsilon
\Delta^\prime_\epsilon$ (to first order in $\epsilon$).
 (d) The deformed counit and antipode of ${\cal{A}}$ are given by
${\cal E}_1=0$ and
$$ S_1(x)=[Ke_{-\rho}e_\sigma,~S(x)]~, \quad x\in {\cal{A}}~.
$$
\no (e) The counit ${\cal E}_\epsilon$ and the antipode $S_\epsilon$
induce analogous structures on ${\cal{A}}/I$.
\vskip.5cm

\no {\bf Proof.}  (a) By the Jacobi identity.  (b) Obvious, for
$\Delta_1(C)=[\Delta(C),~Ke_{-\rho}\otimes Ke_\sigma]\in I\otimes
{\cal{A}}+{\cal{A}}\otimes I$.  (c) Completely straightforward. (d) We
have $({\cal E}\times id)\Delta_1(x)=0$, whence
${\cal E}_1=0$, while
$$ m(id\times S_1)\Delta(H_a) + m(id\times S) \Delta_1(H_a)=0
$$
\no since
$$  m(id\times S_1)\Delta(H_a)=S_1(H_a)= [H_a,~Ke_{-\rho}e_\sigma]~,
$$
$$
\eqalign{m(id\times S)\Delta_1(H_a) &= m(id\times S)(H_a(\sigma)-
H_a(\rho) Ke_{-\rho}\otimes Ke_\sigma \cr &=
\bigl(H_a(\sigma)-H_a(\rho)\bigr) Ke_{-\rho}\bigl(-e_\sigma
e^{-\varphi(\sigma,\cdot)}K^{-1}\bigr) \cr &=
-\bigl(H_a(\sigma)-H_a(\rho)\bigr) Ke_{-\rho}e_\sigma = -[H_a,Ke_{-\rho}
e_\sigma] \cr}
$$
\no and
$$ m(id\times S_1)\Delta(e_\alpha)+m(id\times S)\Delta_1(e_\alpha)=0
$$
\no since
$$
\eqalign{m(id\times S_1)\Delta(e_\alpha) &= S_1(e_\alpha)+e_\alpha
S_1(e^{\varphi(\alpha,\cdot)})\cr &= S_1(e_\alpha)-e_\alpha 
\bigl[ e^{-\varphi(\alpha,\cdot)},~Ke_{-\rho}e_\sigma\bigr] \cr &=
[e_\alpha, Ke_{-\rho}e_\sigma]e^{-\varphi(\alpha,\cdot)}~. \cr m(id
\times S) \Delta_1(e_\alpha) &= -
[e_\alpha,Ke_{-\rho}e_\sigma]e^{-\varphi(\alpha,\cdot)}.\cr}
$$
\no These last two results require some work.
\s (e) This is clear, since ${\cal E}_1=0$ and $S_1(I)\in I$. The proposition 
is proved.
\ve

\no {\bf 15. Exact Deformations of Standard, Generalized Quantum Groups.}

We return to the first order deformations described in Theorem 13.1.
 A deformation of this type,   involving a single pair $(\rho, \sigma)$ for which (13.3) holds,
is called an elementary deformation. We shall see that, to each elementary, first order 
deformation, there is an exact deformation (to all orders in $\epsilon$), that can 
be expressed in closed form. To first order in $\epsilon$, the problem being then
linear, one obtains a more general space of deformations by adding the 
contributions of several  such pairs,
$$
R_1 =  \sum_{(\sigma,\rho) \in [\tau]} 
\,\,\bigl(Rf_\sigma\otimes f_{-\rho} - f_{-\rho}\otimes f_\sigma R -\bigr).
\eqno(15.1)
$$ 
Here the sum is over a subset   $[\tau]$ of the  pairs  $(\sigma, \rho); \,\,\sigma \in \hat
\Gamma_1,\,\,\rho\in\hat\Gamma_2 $, where $\hat \Gamma_{1,2}$ are subsets of the set of 
positive generators, and
$$
 f_\sigma := e^{-\varphi(\sigma,\cdot)}e_\sigma,\,\,\, f_{-\rho} :=
 e_{-\rho}e^{\varphi(\cdot,\rho)},\quad 
e^{\varphi(\sigma,\cdot) + \varphi(\cdot,\rho)} = 1,\quad (\sigma,\rho) \in [\tau].\eqno(15.2)
$$
Not all such compounded, first order deformations
lift to exact deformations.  

The deformed co-product was
also calculated to first order in $\epsilon$, and the results suggest an approach to the 
 exact deformations. The formula (13.4) for $R_1$, as well as the expression (14.4) for  
 the first order deformation of the coproduct,
  both suggest that the deformation be formulated as a twist [D3], but of a type 
much more general than that proposed by Reshetikhin [R]. Additional support for this is found 
in the fact that the exact, elementary deformations 
mentioned above and given below (see ``Examples") are also of this type.
 For the following result ${\cal A}'$ is any coboundary Hopf algebra.
\vskip.5cm

\no {\bf Theorem 15.1.} Let $R$ be the R-matrix, $ \Delta$  the coproduct, of a coboundary
Hopf algebra ${\cal A}',$ and $F \in {\cal A}'\otimes {\cal A}'$, invertible, such that
$$
\bigl((1 \otimes \Delta_{21}) F\bigr)F_{12} = \bigl(
(\Delta_{13} \otimes 1) F\bigr) F_{31}.
\eqno(15.3)
$$
\vskip-.2cm
Then
\vskip-.3cm
$$
\tilde R := (F^t)^{-1} R F\eqno(15.4)
$$
(a) satisfies the Yang-Baxter relation and (b) defines a Hopf algebra $\tilde {\cal A}$ with
the same product and with co-product
$$
\tilde \Delta = (F^t)^{-1} \Delta F^t.\eqno(15.5)
$$

\ve
\no {\bf Proof.} (a) We substitute (15.4) into the expression $\tilde R_{12}\tilde
  R_{13}\tilde R_{3}$. Then  use (15.3) to express $F_{12}(F_{31})^{-1}$ in terms of the
co-products, and the intertwining property of $R$ ($\Delta R = R \Delta'$) to shift the
latter to the  ends. The rest is obvious. (b) It is clear that $\tilde \Delta$ is an
algebra homomorphism. We shall show that the twisted coproduct defined by $\tilde
\Delta$ is co-associative:
$$
\eqalign{
(1 \otimes \tilde \Delta_{23}) \tilde \Delta(x) &= F^{-1}_{32}(1 \otimes \Delta_{23})
\tilde \Delta(x)F_{32}\cr &= F^{-1}_{32}(1 \otimes \Delta_{23} \,(F^t)^{-1})(1 \otimes
\Delta_{23}\, \Delta(x))(1  \otimes \Delta_{23}\, F^t) F_{32},\cr
(\tilde \Delta_{12} +1) \tilde \Delta(x) &= F^{-1}_{21}(\Delta_{12} \otimes 1
\,(F^t)^{-1})(\Delta_{12} \otimes 1 \,\Delta(x))(\Delta_{12} \otimes 1 \,F^t) F_{21}.
\cr}
$$
   Comparing the factors at either end one finds that they agree by virtue of (15.3). 
The result follows, in view of the co-associativity of
$\Delta$.  The theorem is
proved. \footnote*{The connection between Eq.(15.3) and co-associativity was pointed out to 
me by Masaki Kashiwara. The relation makes $F$ a cocycle in the sense of Gerstenhaber [G].}    
\vskip.5cm
\b
 We return to our subject, with $R$ again denoting the standard R-matrix of the algebra
${\cal A}' = {\cal A}/I$. We show first that interesting solutions of (15.3) exist. 
Then we do some preliminary calculations that help us make a general ansatz for $F$ in the
form of a double expansion, $ F = \sum  \epsilon^{nm}F_n^m$, and finally we derive a
recursion relation for $F_n^m$ that will allow us to calculate the classical limit. 
\vskip.5cm
\no {\bf Examples.} An exact deformation of $R$, with the first
order term
$R_1$ as in (15.1) but with the sum reduced to a single term, is given by
$$
F =  e_q^{-\epsilon f_\sigma \otimes f_{-\rho}},\eqno(15.6) 
$$
with
$$
f_\sigma := e^{-\varphi(\sigma,\cdot)}e_\sigma,\,\,\, f_{-\rho} :=
 e_{-\rho}e^{\varphi(\cdot,\rho)}\eqno(15.7)
$$
The q-exponential is   as follows: 
 $q = e^{\varphi(\sigma,\rho)}$, ${e_q}^A := \sum A^n/[n!]_q,\,\, [n!]_q = [1]_q \ldots
[n]_q,\,\, \break [n]_q = (q^n-1)/(q-1)$.  Note that, if $AB = qBA$, then $e_q^A e_q^B =
e_q^{(A+B)}$.   Proposition 15.1 shows that an elementary twist $F$, of the simple
form (15.6), can be combined in a naive way with another elementary twist $\tilde F$, of the same 
type but with $(\sigma,\rho)$ replaced by $(\sigma',\rho')$, only if
$\tilde\Delta(f_\sigma'), \tilde\Delta (\rho')$ reduce to $\Delta(f_\sigma'),
 \Delta (\rho')$; that is, only when the four generators quommute \footnote *{Two elements $X,y$ of an algebra 
quommute if there is $q$ in the field such that $xy - qyx = 0$.} among themselves. 

\vskip.5cm
\no {\bf Notation.} From now on it will be convenient to use the generators $f_{\pm \alpha}$
defined in (15.7). The standard co-product then takes the form
$$
\Delta f_\sigma =  K^\sigma \otimes f_\sigma + f_\sigma \otimes 1, \quad
\Delta f_{-\rho} = 1 \otimes f_{-\rho} + f_{-\rho} \otimes K_\rho,
$$
with
$$
K_\rho := e^{\varphi(\cdot,\rho)},\,\,\,K^\sigma := e^{-\varphi(\sigma,\cdot)}.
$$

The general case of compound deformations is much more complicated. The calculations are
manageable only so long as $F$ can be constructed from  elements of the type $f_\sigma
\otimes f_{-\rho}$ only,  with the factors in this order. A general result is 
Theorem 15.2 below. We need some preparation.
 
\vskip.5cm
 \no {\bf Proposition 15.1.} Let $R_\epsilon$ be an exact deformation of the type   
$$
\eqalign{
&R_\epsilon = (F^t)^{-1}RF, \quad \quad F =  
\sum \epsilon^n  (F_n + \ldots ),\cr
&F_n = \sum_{(\sigma,\rho) \in [\tau]} F_{(\sigma)}^{(\rho')} f_{\sigma_1}\ldots
f_{\sigma_n}  \otimes f_{-\rho_1'} \ldots f_{-\rho_n'},\cr}
\eqno(15.8)
$$
where $+ \ldots $ stands for terms with less than $n$ factors.
  Let
$\Gamma_1, \Gamma_2$ be the subalgebras of ${\cal A}'^+$ generated by $\hat\Gamma_1,
\hat\Gamma_2$. Then we have:  (a) 
There is an isomorphism $\tau: \Gamma_1 \rightarrow \Gamma_2$, such
that the set $[\tau]$ is the restriction of the graph of $\tau$ to $\hat\Gamma_1,
\hat\Gamma_2$,
$$
[\tau] = \{\sigma, \rho \,|\, \sigma \in \hat\Gamma_1, \, \rho = \tau \sigma\in
\hat\Gamma_2\}.\eqno(15.9)
$$
(b) The elements $F_n$ satisfy the recursion relations
$$
 [F_n, f_{-\sigma} \otimes 1] =  ( K^\sigma \otimes f_{-\rho})F_{n-1}
 - F_{n-1}(K_\sigma
\otimes f_{-\rho}),  \quad 
(\sigma,\rho) \in [\tau],\eqno(15.10)
$$
\vskip-2mm
\no 
as well as
\vskip-4mm
$$
 [1 \otimes f_\rho, F_n] = F_{n-1}\bigl( f_\sigma \otimes K^\rho
\bigr) - \bigl(f_\sigma \otimes
K_\rho\bigr) F_{n-1} .\eqno(15.11)
$$
 (c) These recursion relations have the unique solution  
$$
F_{(\sigma)}^{(\rho')} = (-)^n\,\overline t_{(\sigma)}^{(\sigma')},\,\,\,\, (\rho'_1, \ldots ,
\rho'_n) =
\tau(\sigma'_1, \ldots , \sigma'_n),\eqno(15.12)
$$
where the coefficients on the right are the same as in Eq.(2.5), except that $\varphi$ is
replaced by $\varphi^t$. 
\vskip.5cm 
\no {\bf Proof .} We begin by offering some justification for the assumptions.  In view
of the form of $R_1$ it is expected that
$R_n$ is a sum of products of factors of three types:
$$
e_{-\alpha} \otimes e_\alpha,\,\,\,f_{-\rho} \otimes f_{\sigma},\,\,\, 
f_{\sigma} \otimes f_{-\rho},\quad \sigma \in \hat\Gamma_1,\, \rho \in 
\hat\Gamma_2,\eqno(15.13)
$$ 
with coefficients in ${\cal A}' \otimes {\cal A}'$.  In $R_n$
, we isolate the terms with the highest number of factors of the third type, 
$$
X_n = \sum A_{(\sigma)}^{(\rho')}\bigl(e_{-\alpha_1} \ldots e_{-\alpha_k} \otimes
e_{\alpha_1} \ldots e_{\alpha_k}\bigr) B_{(\sigma)}^{(\rho')} \bigl(f_{\sigma_1} \ldots
f_{\sigma_n} \otimes f_{-\rho_1'} \ldots f_{-\rho_n'}\bigr),  
$$
   We shall show that $R_n$contains $X_n \neq 0$.

Let 
$$
Y\hskip-1mm B_\epsilon := R_{\epsilon \,12} R_{\epsilon\,13} R_{\epsilon\,23} - 
R_{\epsilon\,23} R_{\epsilon\,13} R_{\epsilon\,12} \in {\cal A}' \otimes {\cal A}' 
\otimes {\cal A}'. \eqno(15.14)
$$
All terms in $Y\hskip-1mmB_\epsilon$ of order $\epsilon^n$, that have $n$ factors of the
third type in spaces 1,2 are contained in
$$
P_n := F_{n\,12} R_{13} R_{23} - R_{23} R_{13} F_{n\,12}.\eqno(15.15)
$$
For these terms to cancel among themselves $X_n$ must take the form 
$$
X_n = R F_n,\quad F_n = F_{(\sigma)}^{(\rho')} \bigl(f_{\sigma_1} \ldots
f_{\sigma_n} \otimes f_{-\rho_1'} \ldots f_{-\rho_n'}\bigr).\eqno(15.16)
$$
The sum is over all pairs  $(\sigma,\rho) \in [\tau]$ and all permutations $(\rho')$ of
$(\rho)$.
 
Next, the recursion relation (15.10) follows easily from the Yang-Baxter relation (more
precisely from an
examination of terms of low order in space 2), and (15.11) from a similar calculation.
  
We have $F_0 = 1$ and $F_1 = \sum f_\sigma \otimes
f_{-\rho}$. Taking $n=1$ in (15.10) or (15.11) one
gets,  
$$
[f_\alpha, f_{-\beta}] =  \delta^\beta_\alpha
\bigl(e^{\varphi(\cdot,\alpha)} - e^{-\varphi(\alpha,\cdot)}\bigr),\eqno(15.17)
$$
which is confirmed by the definitions in (15.7) and the relation (2.4).
When (15.10) is reduced to a recursion relation for the coefficients, then it turns out
to agree (up to a sign and $\varphi \rightarrow \varphi^t$) with the recursion relations (2.14) for the coefficients 
$t_{(\sigma)}^{(\sigma')}$.  
The integrability of these relations is precisely the statement (a) of  the theorem, as
follows easily from the analysis of these recursion relations in Section 5. Finally,  when (a)
holds, then the  relation (15.11) is equivalent to (15.10).
The proposition  is proved.
\vskip.5cm
After these preliminary explorations we are able to formulate a general result.   
\vskip.5cm
\no {\bf Theorem 15.2.} Let $\Gamma_1, \Gamma_2$ be subalgebras of ${\cal A}'^+$,
generated by subsets $\hat \Gamma_1,\hat \Gamma_2$ of the generators, and $\tau: 
\Gamma_1 \rightarrow \Gamma_2$ an algebra isomorphism.   Let $F \in {\cal A}' \otimes {\cal A}'$ be a formal
series of the form
$$
F = 1 + \sum_{n = 1}^\infty \epsilon^n F_n, \quad F_n = \sum F_{(\sigma)}^{(\rho)} 
f_{\sigma_1}\ldots f_{\sigma_n} \otimes f_{-\rho_1} \ldots f_{-\rho_n}.
\eqno(15.18)
$$
The second sum is here over all $\sigma_i \in \hat \Gamma_1,  \rho_i \in \hat\Gamma_2$. (!)
Note that $F_n$ is   a power series in $\epsilon$.
Suppose that $F$ satisfies (15.3), and that
$$
F_0 = 1,\quad F_1 = -\sum_{\tau^m \sigma = \rho}\epsilon^{m-1}(f_\sigma \otimes f_{-\rho}),\eqno(15.19)
$$
then $F_n$ satisfies
$$
(1 \otimes K_\rho \partial_\rho) F_n + \sum_{\tau^m \sigma = \rho}  \epsilon^m \,\,
[1 \otimes f_\sigma,F_n] + \sum_{\tau^m \sigma = \rho}\epsilon^{m-1}(f_\sigma
\otimes K^\sigma) F_{n-1} = 0.\eqno(15.20)
$$
  With $F_0$ and $F_1$ thus
fixed,  $F_2,F_3,\ldots$ are determined recursively and uniquely. (The operator $K_\rho\partial_\rho$
is the derivation that replaces $f_{-\rho}$ by $K_\rho $.)
\vskip.5cm
\no {\bf Notation.} The sums in (15.19-20), and similar sums to follow, should be understood 
to run over $\sigma \in \hat \Gamma_1$ and over all values of the positive integer $m$ such that $\tau^m\sigma$
is defined; that is, all values of $m$ such that $\tau^{m-1}\sigma \in \hat\Gamma_1$.
\vskip.5cm
\no {\bf Proof.} The exact form (15.19) of
$F_1$  can be inferred directly from the Yang-Baxter relation. That Eq.(15.3) implies   (15.20) is a simple calculation; one
collects all terms that have exactly one generator in the second space.   Let  us verify
that the recursion relation is satisfied for
$ n = 1$ by (15.19). The second term is
$$
-\sum_{\tau^{m'} \sigma'  = \rho'}\epsilon^{m'} 
 \sum_{\tau^m \sigma = \rho}\epsilon^{m-1}f_{\sigma'} \otimes [f_{\sigma}, f_{-\rho'}].
$$
The commutator is
$$
 [f_{\sigma}, f_{-\rho'}] = e^{\varphi(\cdot,\rho')} - e^{-\varphi(\rho',\cdot)}
 = e^{\varphi(\cdot,\tau^{m'}\sigma')} - e^{\varphi(\cdot,\tau^{m'+1}\sigma')}.
$$
The double sum reduces to $\sum_{\tau^m \sigma = \rho}\epsilon^{m-1}f_\sigma
\otimes (K^\sigma - K_\rho)$ and (15.20) reduces to an identity. It remains to prove that 
(15.20) has a unique solution.  Consider first the case that $\hat\Gamma_1 \cap \hat\Gamma_2
$ is empty; then the second term in (15.20) vanishes and the third term
reduces to the term $m=1$. The
recursion relation then reduces to the same form as that which determines the coefficients
of the standard R-matrix, which is known to be integrable. (In this case Proposition 15.1
is the complete solution of the problem, for there are no terms ``$+\ldots$" in (15.8).) In
the general case, when 
$\hat\Gamma_1 \cap \hat\Gamma_2
$ can be non-empty, the second term in (15.20) makes the solution  more difficult, 
but the existence of a solution can still be proved. To do this we  expand $F_n$
as a power series in $\epsilon$, with constant term

$$
F_n^1 = \sum F_{(\sigma)}^{(\tau\sigma')} f_{\sigma_1'}\ldots f_{\sigma_n'} \otimes 
f_{-\tau\sigma_1'}\ldots f_{-\tau \sigma_n'},
$$
and determine the coefficients recursively.
The problem
is therefore  always the integrability of $K_\rho\partial_\rho X = Y, \rho \in \hat\Gamma_2$,
with $Y \in {\cal A}'$  
 given, and this we know  to have a unique solution in ${\cal A}'$, as already noted. 
The theorem is proved.  
\vskip.5cm
The converse, that the solution of (15.20) with $F_0=1$ and $F_1$ given by (15.19) satisfies
(15.3) (and therefore gives a solution of the Yang-Baxter relation) was proved only in the
special case that 
$\hat\Gamma_1 \cap \hat\Gamma_2$ is empty.   Further direct computation 
supports the idea that $R_\epsilon$ always has the form $(F^t)^{-1}RF^t$, with F of the 
form assumed in (15.18). This is strong support for the belief that the solution of the 
recursion relation (15.20), which was proved to exist always, actually furnishes the solution
to the problem of exact deformations in the general case. The results stated in Theorems 13.1 
and 13.2 may 
also be considered as strong evidence.
As we shall see, additional favorable evidence comes from an examination of the classical
limit. To prepare for this we need
\vskip.5cm
\no {\bf Proposition 15.2.} Let 
$$
F_n^m = \sum_{\rho = \tau^m \sigma}  \overline t_{(\sigma)}^{(\sigma')}
f_{\sigma_1}
\ldots f_{\sigma_n}
\otimes f_{-\tau^m \sigma_1'} \ldots f_{-\tau^m \sigma_n'},\quad F^m_0 = 1, 
 \eqno(15.21)
$$
in which the sum extends over $\sigma_i \in \hat \Gamma_1$, $(\sigma')$ a permutation
of $(\sigma)$, and the coefficients $\overline t_{(\sigma)}^{(\sigma')}$ are the same as in
(15.12). Then the unique solution of (15.20) is
$$
\eqalign{
F_n &= \sum_{\Sigma n_i = n} \epsilon^{n_2 + 2n_3 + \ldots} F_{n_1}^1F_{n_2}^2\ldots\ = 
F_n^1 - \epsilon F_{n-1}^1 F^2_1 + \epsilon^2\bigl(F^1_{n-2} F^2_2 + F_{n-1}F^3_1\bigr)
 + \ldots,\cr
F &= \sum \epsilon^nF_n = \sum \epsilon^{n_1 + 2n_2 + \ldots}F^1_{n_1}F^2_{n_2}
\ldots = F^1F^2\ldots,\quad F^m = \sum \epsilon^{nm} F_n^m.\cr}\eqno(15.22)
$$
\vskip1.5cm

\no {\bf 16. Esoteric r-matrices.}
\s
\no \underbar {a) Quantized Kac-Moody algebra
of finite type.}  
\vskip.5cm
\no {\bf Proposition 16.1.} If ${\cal A}'$ is a quantized Kac-Moody algebra of finite type, then $\hat\Gamma_1$
is a proper subset of the set of positive generators and $\tau^{m+1}\hat\Gamma_1 \cap \hat\Gamma_1$ is a 
proper subset of $\tau^m\hat\Gamma_1 \cap \hat\Gamma_1$.
\vskip.5cm
\no {\bf Proof.} Suppose that the statement is false. Then there is $f_\sigma \in \hat\Gamma_1$ such that 
$\tau^m f_\sigma \in \hat\Gamma_1$ for all $m$, and consequently $ \tau^kf_\sigma = f_\sigma$ for some k.
But the condition (15.2) on the parameters, in the classical limit, implies that
$$
\varphi(\tau^m\sigma,\cdot) + \varphi(\cdot,\tau^{m+1}\sigma) = 0.
$$ Summing over $m=0,1,\ldots,k-1$ we obtain
$$
\sum_m(\varphi + \varphi^t)(\tau^m\sigma) = 0,
$$
which contradicts the fact that the Killing form is non-degenerate.
\vskip.5cm
In the classical limit
$$
R_\epsilon = 1 + \hbar r_\epsilon + o(\hbar^2)\quad r_\epsilon = r + \epsilon +
o(\epsilon^2).\eqno(16.1)
$$
In the case of an exact elementary deformation $R_\epsilon$,   the associated  exact
deformation $r_\epsilon$ of $r$ coincides with the first order,
$$
r_\epsilon = r + \epsilon r_1.  \eqno(16.2)
$$

Consider the general case of an exact deformation of $R$ of the form postulated in  
Theorem 15.2. Define $X_\epsilon$  by  
$$
F = 1 + \hbar X_\epsilon + o(\hbar^2), \eqno(16.3)
$$
\vskip-.2cm \no so that \vskip-.5cm 
$$
r_\epsilon = r + X_\epsilon - X^t_{\epsilon}.\eqno(16.4)
$$
\vskip.5cm
\no {\bf Notation.} In this section the symbols $\Gamma_{1,2}$ stand for Lie algebras, 
the classical limits of the algebras so designated until now.
\vskip.5cm

From the fact that the coefficients in the expansion (15.21) of $F$ are the same as the coefficients in the expansion
 (2.5) of the standard R-matrix,
and the known classical limit of the standard R-matrix for a Kac-Moody algebra of finite
type, we get without calculation that
$$
X_\epsilon = -\sum_m \sum_{E_i \in \Gamma_1}\epsilon^{nm} E_i \otimes E_{-\tau^m i},\eqno(16.5)
$$
in which $n$ is the height of $E_i$. The normalization is the same as in
Sections 9-10; more precisely it is fixed as follows. (a) The set $\{E_i\}$ includes the 
generators of $\Gamma_1$. (b) The statement (9.9).
 \footnote*{ Condition (b) can be re-phrased as follows. Let $\Gamma_1^-$ be
the  Lie algebra generated by $\{f_{-\sigma}\},~f_\sigma \in \hat \Gamma_1$ and $\Gamma$ the 
Lie 
algebra generated by $\{f_{\pm \sigma}\},~f_\sigma \in \hat\Gamma_1$. Then \break
$
\sum_{E_i \in \Gamma_1} E_i \otimes E_{-i}
$
is the projection on $\Gamma_1 \otimes \Gamma^-$ of a $\Gamma$-invariant element of 
$\Gamma \otimes \Gamma$.} 
Consequently, 
$$
r_\epsilon = r -  \sum_m \sum_{\matrix {E_i \in \Gamma_1\cr E_j = \tau^m E_i
\cr}} \epsilon^{nm}
\,\, E_i 
\wedge E_{-j}.\eqno(16.6)
$$
The sums are finite, by Proposition 16.1.  
A renormalization exists that reduces the numerical coefficients to 
unity ($\epsilon$  
 in ~\Crm); the result is in complete agreement with [BD].

\ve

\no \underbar {b) Deformations in the affine case.}
Let ${\cal A}'$ be a quantized Kac-Moody algebra of affine type. Two cases should be
distinguished. If the subsets $\hat\Gamma_{1,2}$ of positive roots do no include the 
imaginary root $e_0$, then the formula (16.6) applies without change, except that now $r$ 
is   one of the standard affine r-matrices determined earlier, Eq.s~(9.15), (10.5),
(11.4) or (11.5). There
is nothing more to be said about this case and we turn our attention to the other one.

 What  merits special
attention is the possibility that the first order deformation (15.1)  may include one of
the following
$$
e_0 \wedge e_{-\rho} =  \mu(E_- \otimes e_{-\rho}) -\lambda (e_{-\rho} \otimes E_-),
\eqno(16.7)
$$ 
or
$$ 
 e_\sigma \wedge e_{-0} = \lambda^{-1} (e_\sigma \otimes E_+) - \mu^{-1}(E_+ \otimes
e_\sigma),\eqno(16.8)
$$
with 
$$
\varphi(\cdot,\rho) + \varphi(0,\cdot) = 0, \quad {\rm resp.} \quad \varphi(\cdot,0) +
\varphi(\sigma,\cdot) = 0,\eqno(16.9)
$$
which implies that $\rho  \neq 0$, resp. $ \sigma \neq 0$. A simple renormalization, that
connects the principal picture to the homogeneous picture, brings   (16.8) to the form
$$
e_\sigma \wedge e_{-0} = \sqrt  {\mu/\lambda} (e_\sigma \otimes E_+) -
\sqrt{\lambda/\mu}(E_+ \otimes e_\sigma).
$$ 
  
To deal with the general case of exact deformations it is useful to note the following
\vskip.5cm
\no {\bf Proposition 16.2} If ${\cal A}'$ is a quantized Kac-Moody of affine type, then {\it either }the 
statement about $\hat\Gamma_1$ in Proposition 16.1 continues to hold, {\it or } ${\cal A}'$ is of type 
$A_{N-1}^{(1)}$, $\hat\Gamma_1$ consists of all the positive generators, and $\tau$ generates the cyclic 
group of order $N$.
\vskip.5cm
\no{\bf Proof.} Suppose there is $f_\sigma \in \hat\Gamma_1$ such that $\tau^Nf_\sigma = f_\sigma $ for some $N$. 
Then the Killing form is degenerate. But it is known [K] that any subalgebra of a Kac-Moody algebra of affine
type, obtained by removing one generator, is a Kac-Moody algebra of finite type. It follows that  
  $\hat\Gamma_1$ contains all the positive generators and exactly one $\tau$ orbit. 
Then $\hat\Gamma_1 = \hat\Gamma_2$ and 
$\tau$ lifts to an isomorphism of the Dynkin diagram, which implies the result.
\vskip.5cm  
\no In this section we exclude the exceptional case. This means that the classical limit of $\Gamma_1$
is a finite dimensional Lie algebra, so that (16.6) can be applied directly, since the sum is
finite.

Alternatively, the classical limit   can be found with the 
help of the recursion relation  
 $$
(1 + K_\rho\partial_\rho)F^m_n = -  (f_\sigma \otimes K_\rho)F^m_{n-1},\quad 
\tau^m \sigma = \rho,\eqno(16.10)
$$
or better, the equivalent relation
$$
[1 \otimes f_\rho, F^m_n] =    -  \biggl((f_\sigma \otimes K_\rho) F^m_{n-1} - 
F^m_{n-1}(f_\sigma
\otimes K^\rho)\biggr)\eqno(16.11)
$$
for $F^m_n = \delta_n^0 + \hbar X^m_n + o(\hbar^2)$. This implies that $X^m = 
\sum_{n = 0,1,\ldots} \epsilon^{mn} X^m_n$ (a finite sum) is the
unique solution (of the form that appears in (16.5)) of
$$
[1 \otimes f_\rho +  \epsilon^m f_\sigma \otimes 1, X^m] =  \epsilon^m  
f_\sigma \otimes  (\varphi + \varphi^t)(\rho),   \quad \tau^m\sigma = \rho \in
\hat\Gamma_2.\eqno(16.12)
$$
 \vskip.5cm
\no {\bf Example.} Let ${\cal A}'_{cl}$ be the untwisted, affine Kac-Moody algebra 
$ \widehat{sl(N)}$. A set of positive Serre generators is provided by
the unit matrices $e_i = e_{i,i+1}, ~ i = 1, \ldots N-1$. Set $e_N = e_0 = \lambda e_{N1}$.
The ``most esoteric" deformation (the one with the largest $\Gamma_1$) is defined 
as follows. Take $\Gamma_1$ to be generated by $e_i,~i = 1, \ldots N-1$, 
and $\tau e_i = e_{i+1},~i = 1,\ldots N-1$.  Then $X^m = \sum_n \epsilon^{nm}X_n^m$
with 
$$
X^m_n = -\sum_{i+m+n \leq N}e_{i,i+n} \otimes e_{i+m+n,i+m} -\sum_{i+m+n = N+1} e_{i,i+n} \otimes \lambda^{-1} e_{1,i+m} 
$$ 
and \vskip-.3cm
$$
r_\epsilon = r + \biggl( \sum \epsilon^{nm} X^m_n - 
{\rm transpose}\biggr).
$$
Taking $N=3$ one obtains
$$
r_\epsilon = r - \biggl(\epsilon \,e_{12} \otimes e_{32} + 
\epsilon^2 e_{13} \otimes \lambda^{-1} e_{12} + \epsilon^2 e_{12} \otimes \lambda^{-1}
e_{13} - {\rm transpose}\biggr),
$$
and the renormalization $e_{ij} \rightarrow \lambda^{{j-i \over 3}}e_{ij}$ gives the final 
result
$$
\eqalign{
r_\epsilon = r &- \epsilon \{ \xi^{-1}e_{12} \otimes e_{32} 
+ \xi^{-1} e_{23} \otimes e_{13} + \xi^{-2} e_{13} \otimes e_{12}\}  
- \epsilon^2
\xi^{-1}e_{12} \otimes e_{13}\cr & + \epsilon \{ \xi e_{32} \otimes e_{12} 
+ \xi  e_{13} \otimes e_{23} + \xi^2 e_{12} \otimes e_{13}\} + \epsilon^2
\xi e_{13} \otimes e_{12},\cr }
$$
with $\xi = (\lambda/\mu)^{1/3}$. The un-deformed piece is
$$
r = \varphi + \sum_{i<j}e_{ij} \otimes e_{ji} = {1\over 3}\big(\sum e_{ii} \otimes e_{ii} - 
e_{11} \otimes e_{22} - e_{22} \otimes e_{33} - e_{33} \otimes e_{11}\bigr)+ \sum_{i<j}e_{ij} \otimes e_{ji},
$$
$\varphi$ being   fixed by the relations (15.2). 
This is in agreement with [BD], after transposition and setting $\xi = e^{u/3}, \epsilon = 1$.
.

\ve

\no {\bf 17. Universal Elliptic R and r-matrices.} 

Here we consider the exceptional case (Proposition 16.2) in which $\hat\Gamma_1$
contains all the generators of ${\cal A}'^+$, ${\cal A}'$ is of type $A_{N-1}^{(1)}$ and $\tau^N = 1$.

The expression (15.19) for $F_1$ can be justified as before and the sum is convergent if we interpret $\epsilon$
in \Crm ~and stipulate that
$$
|\epsilon| < 1,
$$
namely
$$
F_1 = {-1\over 1 - \epsilon^N}\sum_{m=1}^N \sum_{\sigma \in \hat\Gamma_1} \epsilon^m 
f_\sigma \otimes f_{-\tau^m\sigma}.\eqno(17.1)
$$
Most, but not all, of the infinite sums that arise can be made meaningful in this way. 
In particular, (15.20) becomes 
$$
(1-\epsilon^N)(1 \otimes K_\rho\partial_\rho)F_n + \sum_{m=1}^N \epsilon^m\,
[1 \otimes f_{\tau^{-m}\rho},F_n] + 
\sum_{m=1}^N \epsilon^{m-1}(f_{\tau^{-m}\rho} \otimes K^{\tau^{-m}\rho})F_{n-1}.\eqno(17.2)
$$
We verify directly that it holds for $n=1$. The second term is
$$\eqalign{
{-1\over 1-\epsilon^N}\sum_{n=1}^N \sum_{m=1}^N&\epsilon^{m+n}f_{\tau^{-m-n}\rho}
\otimes (K_{\tau^{-m}\rho} - K_{\tau^{1-m}\rho})\, \cr
&={-1\over 1-\epsilon^N}\sum_{M=1}^N \epsilon^M f_{\tau^{-M}\rho}
\otimes (K^{\tau^{-M}\rho} - K_\rho)(1-\epsilon^N).
\cr}
$$
  The terms $K^{\tau^{-M}},K_\rho$ comes from the ends of the 
summation   while all the other terms cancel pairwise since $K^\sigma = K_{\tau\sigma}$.

The infinite product
$$
F =  F^1F^2\ldots\eqno(17.3)
$$
cannot be given anything more than a formal significance in the structural context but, as will be shown below, in a finite 
dimensional representation the question of convergence (with $\epsilon $ in \Crm) is not 
difficult. We define $F^m$ by the (always uniquely integrable) relation (16.10),
$$
(1 \otimes K_\rho\partial_\rho)F^m = -\epsilon^m(f_{\tau^{-m}\rho}\otimes K_\rho)F^m,\quad F^m = 1 - \epsilon^m
\sum_\sigma f_\sigma \otimes f_{-\tau^m\sigma} + o(\epsilon^{2m}),\eqno(17.4)
$$
or its equivalent
$$
[1 \otimes f_\sigma,F^m] = -\epsilon^m\biggl((f_{\tau^{-m}\sigma} \otimes K_\sigma)F^m
-F^m(f_{\tau^{-m}\sigma} \otimes K^\sigma)\biggr),\eqno(17.5)
 $$  
with the same initial condition. We verify that, with this definition of $F^m$, (17.3) satisfies (17.2) or
$$
(1 - \epsilon^N)(1 \otimes K_\rho\partial_\rho)F + \sum_{\tau^n\sigma=\rho} \epsilon^n[1 \otimes f_\sigma,F]
+\sum_{\tau^n\sigma=\rho} \epsilon^n(f_\sigma \otimes K^\sigma)F = 0.\eqno(17.6)
$$
The range of the summation is $ n =  1,2,\ldots,N,\, \sigma \in \hat\Gamma_1$.  One has
$$
\eqalign{
&\sum_n \epsilon^n[1 \otimes f_{-\tau^{-n}\rho},F^mF^{m+1}] = 
-\sum_n \epsilon^{m+n} f_{\tau^{-m-n}\rho} \otimes
K_{\tau^{-n}\rho}F^mF^{m+1}\cr
&\hskip.2cm + F^m\biggl\{ \sum_n \epsilon^{m+n} f_{\tau^{-m-n}\rho} \otimes K^{\tau^{-n}\rho} - 
\sum_n \epsilon ^{m+n+1} f_{\tau^{-m-n-1}\rho} \otimes K^{\tau^{-n-1}\rho}\biggr\} F^{m+1} + \ldots \,.\cr}
$$
In the second line everything cancels except for the first and the last terms, leaving
$$
-\sum_n \epsilon^{m+n}( f_{\tau^{-m-n}\rho} \otimes K_{\tau^{-n}\rho}) F^m F^{m+1} 
+ (1-\epsilon^N)F^m \epsilon^{m+1}( f_{\tau^{-m-1}\rho} \otimes K_\rho) F^{m+1} + \ldots\,.
$$
The total contribution of the commutator in (17.6) is thus
$$
-\sum_{n=1}^N \epsilon^{n+1}( f_{\tau^{-n-1}\rho} \otimes K_{\tau^{-n}\rho})F 
+ (1-\epsilon^N)\sum_{m=1}^\infty F^1\ldots F^m \epsilon ^{m+1}(f_{\tau^{-m-1}\rho} \otimes K_\rho) F^{m+1} \ldots \,.
$$
Adding the first term in (17.6) leaves us with
$$
-\sum_n \epsilon^{n+1}( f_{-\tau^{-n-1}\rho}\otimes K_{\tau^{-n}\rho})F - 
\epsilon(1-\epsilon^N) (f_{-\tau^{-1}\rho} \otimes K_\rho) F 
= -\sum_n \epsilon^n (f_{-\tau^{-n}\rho} \otimes K^{\tau^n\rho}) F,
$$
 which is cancelled by the last term.

In the classical limit $F^m = 1 +  \hbar X^m + o(\hbar^2)$ and $X^m$ satisfies (16.12).
We shall solve these relations in the case of the simplest affine Kac-Moody algebra. Set

$$[f_1,f_{-1}] = (\varphi + \varphi^t)(1) = \sigma_3,\quad \tau f_1 = f_0,~~\tau f_0 = -f_1,
$$
and
$$X^m = A^m \sigma_3 \otimes \sigma_3 + B^m(f_1 \otimes f_{-1} + f_0 \otimes f_{-0}) + C^m(f_1 \otimes f_{-0} + f_0 \otimes
f_{-1})
$$
and impose (16.12). The result is, with $x =  \sqrt{\lambda/\mu}$,
$$\eqalign{
A^m &= \sum_{n=1}^\infty (-\epsilon^{2n})^m x^{-n},\cr
B^{2m} &= \sum_{n=1}^\infty(\epsilon^{2n-1})^{2m}x^{1-n}, \quad B^{2m-1} = 0,\cr
C^{2m-1} &= \sum_{n=1}^\infty(\epsilon^{2n-1})^{2m-1}x^{1-n}, \quad C^{2m} = 0.\cr}
$$
The deformed r-matrix is $r_\epsilon = r + X - X^t$, with 
$$\eqalign{
X = \sum_{n=1}^\infty X^m = 
\sum_n {-\epsilon^{2n}\over 1 + \epsilon^{2n}}x^{-n} \sigma_3\otimes \sigma_3  
&+ \sum_{n=1}^\infty {\epsilon^{4n-2}\over 1- \epsilon^{4n-2}}x^{1-n}(f_1 \otimes f_{-1} + 
f_0 \otimes f_{-0})\cr 
&+ \sum_{n=1}^\infty {\epsilon^{2n-1}\over 1- \epsilon^{4n-2}}x^{1-n}(f_1 \otimes f_{-0} + f_0 \otimes f_{-1}).\cr}
$$
Setting $\lambda/\mu = e^{2\pi iu}$ one gets
$$\eqalign{
(i/2)(X-X^t) = &\sum_{n=1}^\infty\biggl\{ {-\epsilon^{2n}\over 1 + \epsilon^{2n}}(\sigma_3 \otimes \sigma_3) \sin 2n\pi u \cr
&\hskip.2cm + {\epsilon^{4n-2}\over 1 - \epsilon^{4n-2}}\bigl(  x \,f_1 \otimes f_{-1} + {1\over x}\, f_{-1} \otimes f_1\bigr) \sin (2n-1)\pi u\cr 
&\hskip.5cm + {\epsilon^{2n-1}\over 1 - \epsilon^{4n-2}}\bigl(\sqrt{1/\mu\lambda}\,f_1 \otimes f_1 + \sqrt{\mu\lambda}\,f_{-1}
\otimes f_{-1}\bigr) \sin (2n-1) \pi u\biggr\}.\cr}
$$
The trigonometric r-matrix (9.15) is
$$
{i\over 2}\biggl({1\over\tan \pi u}(\sigma_3 \otimes \sigma_3) + {1\over   \sin \pi u} \bigl(\sqrt x\,f_1 \otimes f_{-1}
+ \sqrt{1/x}\,f_{-1} \otimes f_1\bigr)\biggr).
$$
Adding, one finds the series expansion of elliptic functions, and complete agreement with the elliptic 
r-matrices of [BD]. To transform to their notation replace 
$$f_1 \rightarrow \sqrt\lambda\, e_{12}, \,f_{-1} \rightarrow
\sqrt{1/\lambda}\, e_{21}\eqno(17.7)
$$

Finally, we shall show that the expression for the Universal Elliptic R-matrix $R_\epsilon = (F^t)^{-1}RF, 
\,F = F^1F^2\ldots$ in terms of  an infinite product is both 
meaningful and useable, by projecting on a finite dimensional represention. We limit ourselves to the 
fundamental representation of $sl(2)$. After  rescaling of the generators as in (17.7), 
$F^m$ and $R_\epsilon$ take the form
$$
F^m = \pmatrix{a^m&&&d^m\cr&b^m&c^m&&\cr
&c^m&b^m&&\cr
d^m&&&a^m},\quad R_\epsilon =  \pmatrix{a&&&d\cr&b&c&&\cr
&c&b&&\cr
d&&&a \cr }.
$$ 
The matrix elements are completely determined by the recursion relation (17.5); namely for  $ m = 1,2,\ldots $,
$$
\eqalign{
a^{2m-1} &= 1 - \epsilon^{4m-2},~~  b^{2m-1} = 1 - \epsilon^{4m-2}{q^2\over x},
 \,\, c^{2m-1} = 0, ~~d^{2m-1} = \epsilon^{2m-1}({1\over q} - q)\sqrt{1\over x},\cr
a^{2m } &= 1 - \epsilon^{4m}{q^2\over x}, \,~~~~~ b^{2m} = 1 - \epsilon^{4m}{1\over x},
 ~~~~~~~\,  c^{2m} = \epsilon^{2m}\sqrt{1\over x}({1\over q} - q),~~~ ~~ d^{2m} = 0,\cr}
$$
and
$$
a+d ~: ~a-d~ :~ b+c ~: ~b-c  = {dn(u+\rho)\over dn(u-\rho) } : 1 : {cn(u+\rho)\over cn(u-\rho)} : {sn(u+\rho)\over sn(u-\rho)}.
$$
A modular transformation brings this into perfect agreement with Baxter [B].

\vskip3.5cm

{\ce {\bf Acknowledgements}}

I am indebted to Moshe Flato for a detailed and constructive criticism  of the original
manuscript, and  Masaki Kashiwara and Tetsuji Miwa for fruitful discussions. 
This work was completed at the Research Institute for Mathematical Sciences,
Kyoto University. I thank the director, Professor Huzihiro Araki, for making the pleasant and productive
stay in Kyoto possible, and for friendly hospitality during my visit. I thank the staff of the
International Office of the Institute for many services, and the Japanese Ministry of
Education, Science, Sports and Culture for financial support. 
\ve

\parindent0pt {\bf References.}

\noindent
\parindent0pt 

[B]\hskip5mm R.J. Baxter, Partition Function of the Reught-Vertex Model, 
\line {\hfill Annals of Physics {\bf 70}  (1972) 193-228.}

[BD]\hskip2mm  A.A. Belavin and V.G. Drinfeld, Triangle Equation and Simple Lie Algebras, 
Sov. \line {\hfil Sci.Rev.Math.{\bf 4} (1984) 93-165.}

[BP]~~ F. Bidegain and G. Pinczon, Quantization of Poisson Lie groups and
applications, 
\line {\hfil Dijon preprint 1995.}

[BFGP]  P. Bonneau, M. Flato, M. Gerstenhaber and G. Pinczon, The hidden group  \line 
{\hfil structure of quantum groups: strong duality, rigidity, and preferred deformations,} 
\line {\hfil Commun.
Math.Phys. (1994) 125-156.}

[C]~~~~  C. Chevalley, ``Sur la classification des algebres de Lie simples
et de leur represen- \line {\hfil tations".    C.R. {\bf 227} (1948)
1136-1138.}

[CG] \hskip1mm  E. Cremmer and J.-L. Gervais, Commun.Math.Phys. {\bf 134} (1990)

[D1]\hskip3mm  V.G. Drinfeld, Quantum Groups, in Proceedings of the International Congress of
   \line {\hfil Mathematicians,~Berkeley,  A.M. Gleason, ed. (American
Mathematical Society,} \line {\hfil Providence  1987.)}

[D2]\hskip3mm  V.G. Drinfeld, in ``Quantum Groups," P.P. Kulish, ed. Proceedings of
the Work- \line { \hfil shop held in the Euler International  Mathematical
Institute, Leningrad, Fall 1990.}

[EK]~~\hskip.8mm P. Etingof and D. Kazhdan, Quantization of Poisson algebraic groups 
and Poisson \line {\hfil homogeneous spaces, Harvard preprint Oct 20 1955.}

[FR]\hskip 3.5mm I.B. Frenkel and N.Yu. Reshetikhin, Quantum Affine Algebras and Holonomic 
   \line {\hfil Difference Equations,Commun.Math.Phys. (1992) {\bf 146} 1-60).}

[F]~~ \hskip1.mm  C. Fr\o nsdal, Cohomology and Quantum Groups, in Proceedings of the 
 Karpacz \line {\hfil  Winter School, February 1994.}

[FG1] C. Fr\o nsdal and A. Galindo, Contemporary Mathematics, {\bf 175}
(1994) 73-88.

[FG2] C. Fr\o nsdal and A. Galindo, Lett.Math.Phys. {\bf 34} (1995) 25-36.

[G] M. Gerstenhaber,

[GK]\hskip2mm  O. Gabber and V.G.  Kac, On defining relations of certain infinite-dimensional Lie
 \line {\hfil   algebras,
Bull. Amer. Math. Soc. {\bf 5} (1980), 185-189. }

  [J] ~~~M. Jimbo, Commun.Math.Phys. {\bf 102} (1986) 537-547.

[K] ~\hskip2mm  V.G. Kac, ``Infinite Dimensional Lie Algebras," Cambridge
University Press 1990.
\ve
[KR]\hskip1.8mm  A.N. Kirillov and N.Yu. Reshetikhin, q-Weyl group and a multiplicative 
formula  \line {\hfil for the universal R-matrix, Commun.Math.Phys.  {\bf 134}
(1991) 421-431.}

[LS]~\hskip2mm S.Z. Levendorski and Y.S. Soibelman, Some applications of quantum Weyl group I.
\line {\hfil The multiplicative formula for universal R-matrix for simple Lie algebras.  }
\line {\hfil J.Geom.Phys. {\bf 7} (1991) 1-14.}

[Ma]\hskip1.8mm  Yu.I. Manin, ``Topics in noncommutative geometry," Princeton
University Press, \line {\hfil Princeton NJ 1991.}

[Mo] \hskip1.4mm R.V. Moody, ``A new class of Lie algebras." J. Algebra {\bf 10}
(1968) 211-230.

[R]~~ \hskip1.35mm  N.Yu. Reshetikhin, Lett.Math.Phys. {\bf 20} (1990) 331-336.

[Ro]~\hskip1.8mm  M. Rosso, An analogue of P.B.W. theorem and the universal R-matrix for $U_hsl(N)$,
\line { \hfil Commun.Math.Phys.{\bf 124} (1989) 307-319.}

[Sc] \hskip2.5mm A. Schirrmacher, Z. Phys.C {\bf 50} (1991) 321.

[Su]~\hskip2mm  A. Sudbery, J.Phys.A.Math.Gen. {\bf 23} (1990) L697.

[T]~\hskip3.5mm T. Tanisaki, Int.J.Mod.Phys.A, Suppl. 1B (1992) 941-961. 

[TV]      P. Truino and V.S. Varadarajan, Quantization of Reductive Lie Algebras, Rev.
  \line {\hfill  Math.Phys. {\bf 5} (1993) 303; 
Universal Deformations of Reductive Lie Algebras,}
 \line {\hfill Lett.Math.Phys.{\bf 26} (1992) 53.}

\end

\vskip 1.5cm

\end